\documentclass[12pt,a4paper]{article}
\usepackage[utf8]{inputenc}
\usepackage[T1]{fontenc}
\usepackage[english]{babel}
\usepackage{amsmath}
\usepackage{amssymb,amstext}
\usepackage{mathrsfs}
\usepackage{amsthm}

\usepackage{multicol}
\usepackage{enumerate}
\usepackage{fancyhdr}

\usepackage{color}
\usepackage{graphicx}
\usepackage{caption}

\usepackage{geometry}

\newtheorem{theorem}{Theorem}
\newtheorem{lem}{Lemma}
\newtheorem{rem}{Remark}
\newtheorem{defin}{Definition}
\newtheorem{prop}{Proposition}
\newtheorem{cor}{Corollary}

\newenvironment{proof1}{\textit{Proof: }}{\begin{flushright}$\blacksquare$\end{flushright}}

\def\bit{\begin{itemize}}
\def\eit{\end{itemize}}
\def\bdes{\begin{description}}
\def\edes{\end{description}}

\def\iti{\item[(i)]}
\def\itii{\item[(ii)]}

\def\itiii{\item[(iii)]}

\def\beq{\begin{equation}}
\def\eeq{\end{equation}}

\def\ben{\begin{enumerate}}
\def\een{\end{enumerate}}

\def\beqar{\begin{eqnarray}}
\def\eeqar{\end{eqnarray}}

\def\beqarr{\begin{eqnarray*}}
\def\eeqarr{\end{eqnarray*}}
\def\RR{{\mathbb R}}  

\def\rar{\rightarrow}
\def\eps{\varepsilon}

\begin{document}
	\title{\textbf{Self-repelling diffusions on a Riemannian manifold}}
	\author{ Michel Bena\"im and Carl-Erik Gauthier\\Institut de Math\'ematiques\\Universit\'e de Neuch\^atel, Switzerland
		\thanks{We acknowledge financial support from the  Swiss National Science Foundation Grant  200020\- 149871/1. We thank B. Colbois and H. Donnelly for useful discussions on eigenfunctions of the Laplace operator and P.Monmarch\'e for useful discussion about hypocoercivity. }}
	\maketitle
	\date{}

	\begin{abstract}
		Let M be a compact connected oriented Riemannian  manifold. The purpose of this paper is to investigate the long time behavior of a degenerate stochastic differential equation  on the state space $M\times \mathbb{R}^{n}$; which is obtained  via a natural change of variable from a self-repelling diffusion taking the form  $$dX_{t}= \sigma dB_{t}(X_t) -\int_{0}^{t}\nabla V_{X_s}(X_{t})dsdt,\qquad X_{0}=x$$ where $\{B_t\}$ is a Brownian vector field on $M$, $\sigma >0$ and $V_x(y) = V(x,y)$ is a diagonal Mercer kernel.
		We prove that the induced semi-group enjoys the strong Feller property  and has a unique invariant probability $\mu$ given as the product of the normalized Riemannian measure on M and a Gaussian measure on $\mathbb{R}^{n}$. We then prove an exponential decay to this invariant probability in $L^{2}(\mu)$ and in total variation.
	\end{abstract}
	\textbf{Keywords:} self-interacting diffusions, strong Feller property, degenerate diffusions, hypocoercivity, invariant probability measure
	
	MSC primary: 58J65, 60K35, 60H10, 60J60, 
	secondary 37A25, 37A30
\section{Introduction}
Let $M$ be a smooth (i.e $\mathcal{C}^{\infty}$) Riemannian manifold,  $V : M \times M \rightarrow \mathbb{R}$ a smooth function and
$w : [0,\infty[ \rightarrow [0, \infty[$ a continuous function.
 Adopting the terminology now coined in the literature   we define a {\em Self Interacting Diffusion with potential $V$ and weight function $w$}
 to be a continuous time stochastic process $(X_t)_{t \geq 0}$ living on $M$  defined by the stochastic differential equation
\begin{equation}\label{eq-1}
dX_{t}= \sigma dB_{t}(X_{t})-  \nabla V_t(X_t) dt,
\end{equation}
where $\sigma>0$, $\{B_{t}\}$ is a Brownian vector field on $M$ and
\begin{equation}\label{eq-1pot}
 V_t(x) = w_t \int_{0}^{t} V(X_s, x)ds,
\end{equation}
The case $M$  compact and $\displaystyle w_t = t^{-1}$  has been thoroughly analyzed in a series of papers by the first named author
in collaboration with Raimond (\cite{BRL}, \cite{BR2}, \cite{BR3}) and Ledoux  \cite{BRL}.
In particular, it was shown that  long term behavior of the normalized occupation measure $\mu_{t} =\frac{1}{t}\int_{0}^{t}\delta_{X_{s}}ds$ can be precisely related to the long term behavior of a deterministic semi-flow defined on the space of probability measures over $M.$ Pemantle's survey paper (\cite{Pem}) contains a comprehensive discussion of these results among others and further references.  Some extensions to noncompact spaces have been considered by Kurtzmann in  \cite{KurtzK}, \cite{Kurtz} and  other weight functions decreasing to zero  by Raimond in  \cite{Rai2}.

When $w$ doesn't converge to zero, say $w_t = 1$, the literature on the subject mainly consists of case studies under the assumption that $M = \RR$ (or $\RR^d$) and $V(x,y) = v(y-x).$
Self attracting processes, that is $ x v'(x) \geq 0$ (or $\langle x, v'(x) \rangle \geq 0$ in $\RR^d$), have been considered by  Cranston and  Le Jan  \cite{CranYLJ}, Raimond \cite{Rai}, Herrmann and Roynette \cite{HR},  Herrmann and Scheutzow \cite{HS} and typically converge almost surely.
For self repelling processes, that is $ x v'(x) \leq 0,$ the process tends to be "transient" and strong law of large numbers and rate of escapes have been obtained under various assumptions by Cranston and Mountford \cite{MountCran}, Durrett and Rogers \cite{DurRog}, Mountford and Tarr\`es \cite{MountTar}. In \cite{TTV12}, Tarr\`es, T\'oth and Valk\'o consider the situation when $v$ is a sufficiently smooth function having a nonnegative Fourier transform. Under this condition and other technical assumptions,
 they show that the {\em environment seen from $X_t$}, that is the mapping
$x \mapsto \int_0^t v'(x + X_t - X_s) ds$, admits  an ergodic invariant Gaussian measure.


 In this paper we will pursue this line of research and investigate the long term behavior of  (\ref{eq-1}) under the assumptions that:
 \bdes
 \iti ({\bf Strong interaction}) $w_t = 1$.
 \itii ({\bf Compactness}) $M$ is smooth, finite dimensional, compact, oriented, connected and without boundary.
 \itiii ({\bf Self repulsion}) 
 $V$ is a\textit{ Mercer kernel}. That is, $V(x,y) = V(y,x)$ and
$$\int_{M}\int_{M}V(x,y) f(x) f(y) dx dy\geqslant 0 $$ for all $f\in L^{2}(dx)$, where $dx$ stands for the Riemannian measure.
\edes
By Mercer Theorem, $V$ can be written as
\beq
\label{eq:mercer}
V(x,y) = \sum_i a_i e_i(x) e_i(y)
\eeq
where $a_i \geq 0$ and $\{e_i\}$ is an orthonormal (in $L^2(dx)$) family of eigenfunctions of the operator $f \mapsto Vf,$ where $Vf(x) = \int V(x,y)f(y) dy.$

Thus, if one interpret the sequence  $$\Psi(x) = (\sqrt{a_i} e_i(x))_i$$ as a {\em feature vector} representing  $x$ in  $l^2$,
$$V(x,y) = \langle \Psi(x), \Psi(y) \rangle_{l^2}$$ can be thought of as a  similarity between the feature vectors  $\Psi(x)$ and $\Psi(y).$
The process is therefore {\em self-repelling} in the sense that the drift term $-\nabla V_t(X_t)$ in equation (\ref{eq-1})
 tends to minimize the similarity between  the current feature vector $\Psi(X_t)$  and the cumulative feature $\int_0^t \Psi(X_s) ds.$


Here we will focus on the particular situation where
\bdes
\item {\bf (iii')}  ({\bf Diagonal decomposition})  The sum in (\ref{eq:mercer}) is finite and the $\{e_i\}$ are eigenfunctions of the Laplace operator.
 \edes
Our  motivation for such a restriction is twofold. First, for a suitable choice of $n$ and $(a_i),$ the feature map
 $$\Psi : M \mapsto \RR^n,$$ $$x \mapsto (\sqrt{a_1} e_1(x), \ldots, \sqrt{a_n} e_n(x))$$ is a quasi-isometric embedding of $M$ in $\RR^n.$
 We refer the reader to the recent paper  (Portegies 2015 \cite{Portegies} ) for a precise statement (Theorem 5.1), and
  further interesting discussions and references on embedding by eigenfunctions.   In particular, for some $\eps > 0$
  $$ - V(x,y) \leq \frac{1}{2} \|\Psi(x) - \Psi(y)\|^2 \leq (1 + \eps) \frac{d(x,y)^2}{2},$$ where $d$ stands for the Riemannian distance on $M.$
  Hence, with this choice of $(a_i),$ the smaller is $V_t(X_t)$ the larger is the cumulative quadratic distance
  $\int_0^t d^2(X_t,X_s) ds.$

  Secondly, under hypothesis $(iii)'$, an invariant probability measure of the process $(X_t, V_t(x))$ can be explicitly computed.
It turns out that this will be of fundamental importance for our analysis.

 \subsection*{\textbf{A motivating example: the periodic case.}}
 Let $M=\mathbb{S}^{1} = \mathbb{R}/2\pi\mathbb{Z}$  denote the unit circle and let $V:M\times M\rightarrow\mathbb{R}$ be the map defined
 by $$V(x,y) =  \cos(y-x) = 1 -\frac{1}{2} d^{2}(y,x),$$ where $d(y,x) = \vert e^{iy}-e^{ix}\vert.$

 Noting that  $\nabla V_{x}(y)=-\sin(y-x),$ (\ref{eq-1}) can be rewritten as
 \begin{equation}\label{eqMotiv}
 dX_{t}=\sigma dB_{t}+\int_{0}^{t}\sin(X_{t})\cos(X_{s})-\cos(X_{t})\sin(X_{s})dsdt.
 \end{equation}
 Setting $U_{t}=\int_{0}^{t}\cos(X_{s})ds$ and $V_{t}=\int_{0}^{t}\sin(X_{s})ds$ we get the following SDE on $\mathbb{S}^{1}\times\mathbb{R}^{2}$:
\begin{equation}\label{eq00}
\left\{
\begin{array}{r c l}
dX_{t}&=& \sigma dB_{t}+(\sin(X_{t})U_{t}-\cos(X_{t})V_{t})dt\\
dU_{t}&=& \cos(X_{t})dt. \\
dV_{t}&=& \sin(X_{t})dt
\end{array}
\right.
\end{equation}
This system enjoys the following properties, summarized by the next Theorem, which proof follows from Theorems \ref{t2},\ref{t4},\ref{t5} and Proposition \ref{p0}.
Given $y = (x,u,v) \in \mathbb{S}^{1}\times\mathbb{R}^{2}$, we let $(Y_t^y)_{t \geq 0} = ((X^y_t,U^y_t,V^y_t))_{t \geq 0}$ denote the solution to (\ref{eq00}) with initial condition $Y_0^y = y.$ Here $\mathbb{S}^{1}$ is identified with $\mathbb{R}/2\pi\mathbb{Z}$.
\begin{theorem}
 The Markov process induced by (\ref{eq00}) is a positive Harris process and admits a unique invariant probability  given as $$\mu(dxdudv)=\frac{dx}{2\pi}\otimes \frac{\exp(-u^{2}/2)}{\sqrt{2\pi}}du\otimes  \frac{\exp(-v^{2}/2)}{\sqrt{2\pi}}dv. $$
Furthermore, the law of $ Y_t^y$ converges  exponentially fast  to $\mu$ in $L^{2}(\mu)$ and in total variation.
\end{theorem}
\begin{rem} A similar result holds for the decoupled SDE when $V(x,y)=\sum_{j=1}^{n}a_{j}\cos (j(y-x))$ and $a_{j}>0$ for all $j=1,\cdots ,n$, by setting $U_{j}(t)=\int_{0}^{t}\cos(jX_{s})ds$ and $V_{j}(t)=\int_{0}^{t}\sin(jX_{s})ds$.
\end{rem}
\begin{theorem} Almost surely, the solution of (\ref{eqMotiv}) with initial condition $(X_{0},U_{0},V_{0})=(0,0,0)$ does not converge on $\mathbb{S}^{1}$ and a fortiori on $\mathbb{R}$. However, on $\mathbb{R}$, $$\frac{X_{t}}{t}\rightarrow 0 \text{ a.s.} \text{ as } t\rightarrow\infty.$$
\end{theorem}
\begin{proof1} Let $\varepsilon >0$ and set $R_{j}^{\varepsilon}=\bigcup_{k\in\mathbb{Z}}((2k+j)\pi-\varepsilon, (2k+j)\pi+\varepsilon)\times\mathbb{R}^{2}$, $j=0,1$. Then by positive Harris recurrence of $(X_{t},U_{t},V_{t})_{t}$, we have that $$X_{t}\in \bigcup_{k\in\mathbb{Z}}((2k+j)\pi-\varepsilon, (2k+j)\pi+\varepsilon),$$
infinitely often for $j=0,1$. This proves the first assertion.\\
Applying now Corollary \ref{corErg} in section \ref{result} to the function $f(x,u,v)=\sin(x)u-\cos(x)v$ gives us
\begin{equation*}
\lim_{t\rightarrow\infty}\frac{1}{t}\int_{0}^{t}f(X_{s},U_{s},V_{s})ds=\int_{\mathbb{S}^{1}\times\mathbb{R}^{2}}f(x,u,v)\mu(dx,du,dv)=0\quad \mathbb{P}_{(0,0,0)} a.s.
\end{equation*}
Consequently,
\begin{equation*}
\frac{X_{t}}{t}=\sigma\frac{B_{t}}{t}+\frac{1}{t}\int_{0}^{t}f(X_{s},U_{s},V_{s})ds
\end{equation*}
converges $\mathbb{P}_{(0,0,0)}$ almost surely to $0$.
\end{proof1}

\begin{center}
\parbox{0.3\linewidth}{\includegraphics[width=0.35\textwidth]{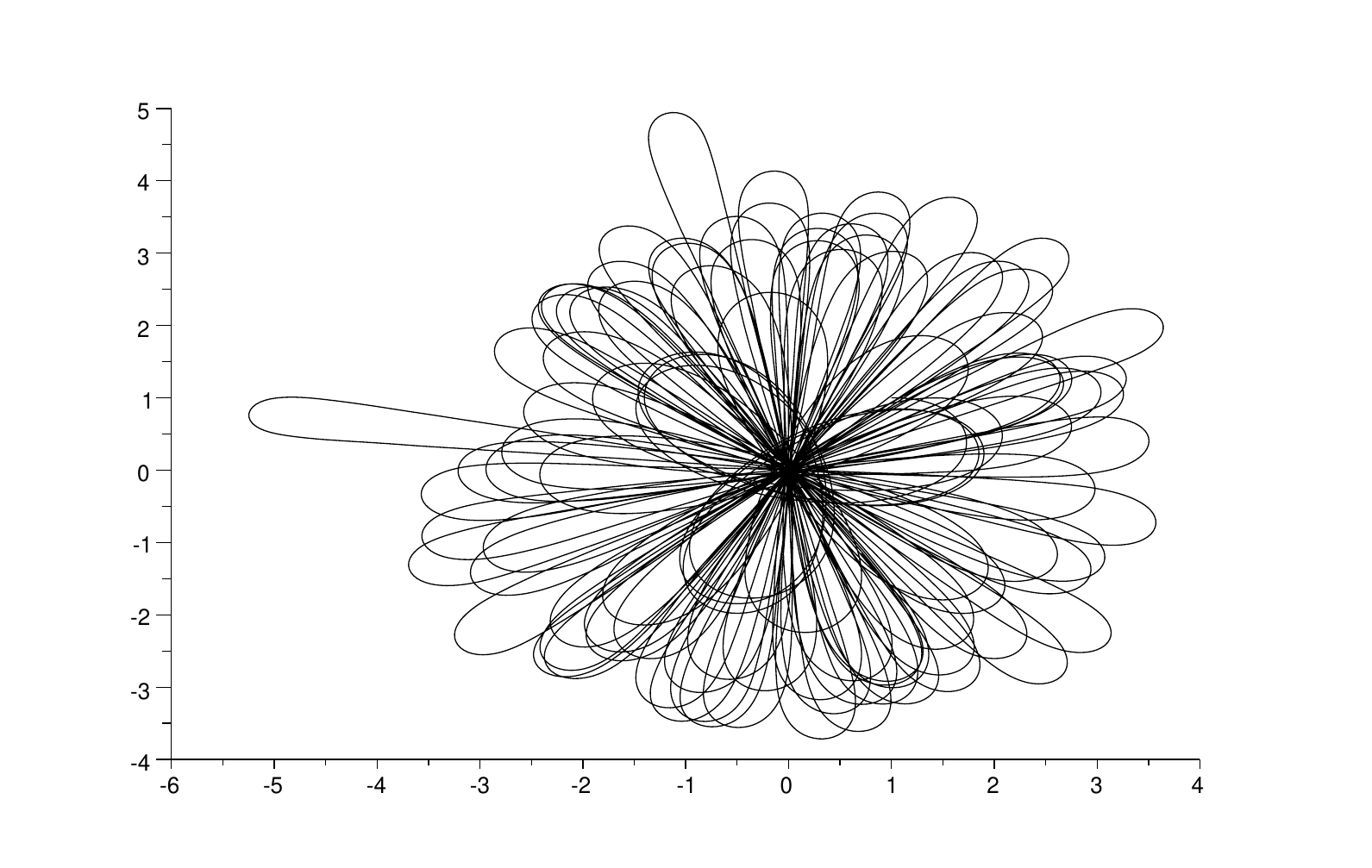}}
\parbox{0.2\linewidth}{}
\parbox{0.3\linewidth}{\includegraphics[width=0.35\textwidth]{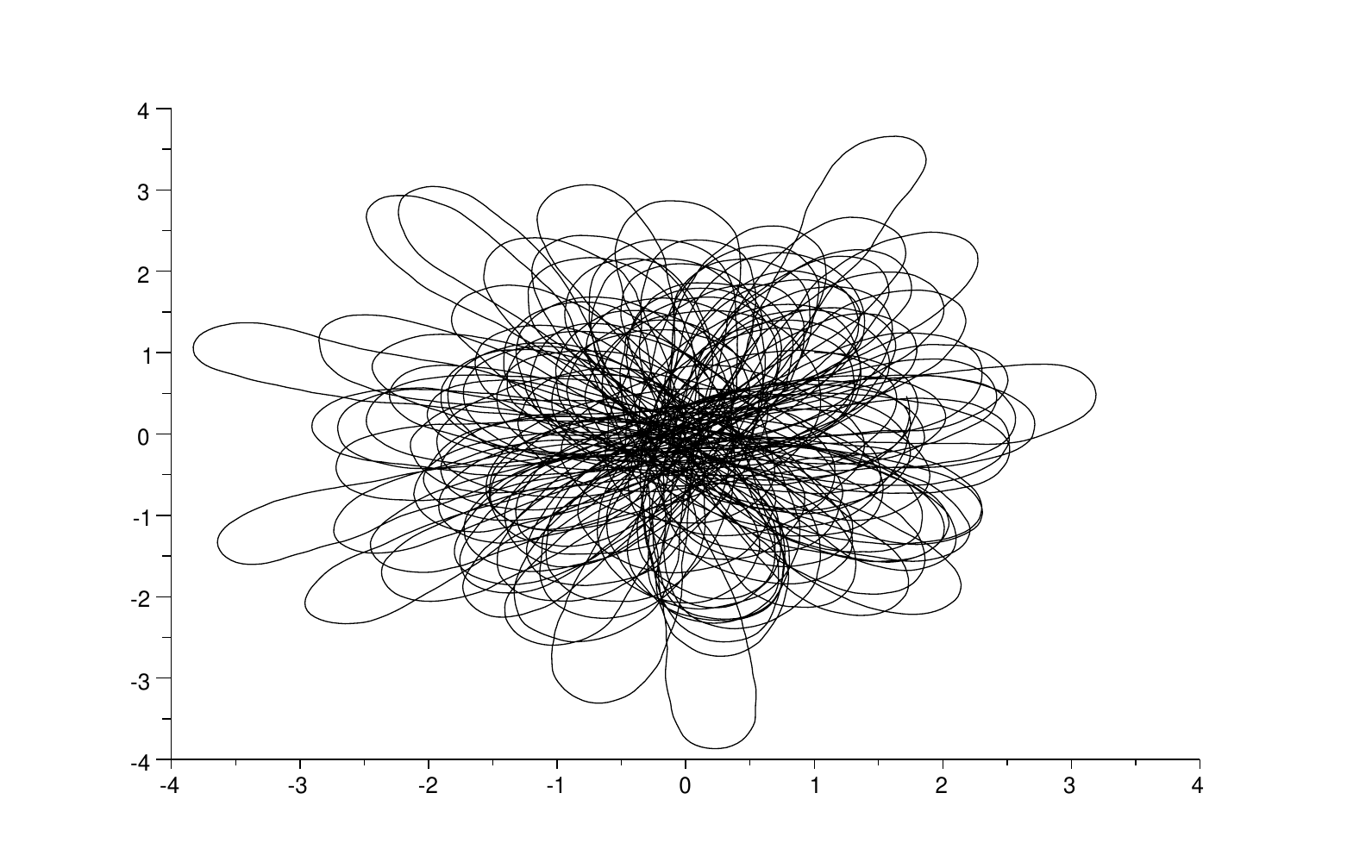}}
\parbox{0.2\linewidth}{}
\parbox{0.3\linewidth}{\includegraphics[width=0.35\textwidth]{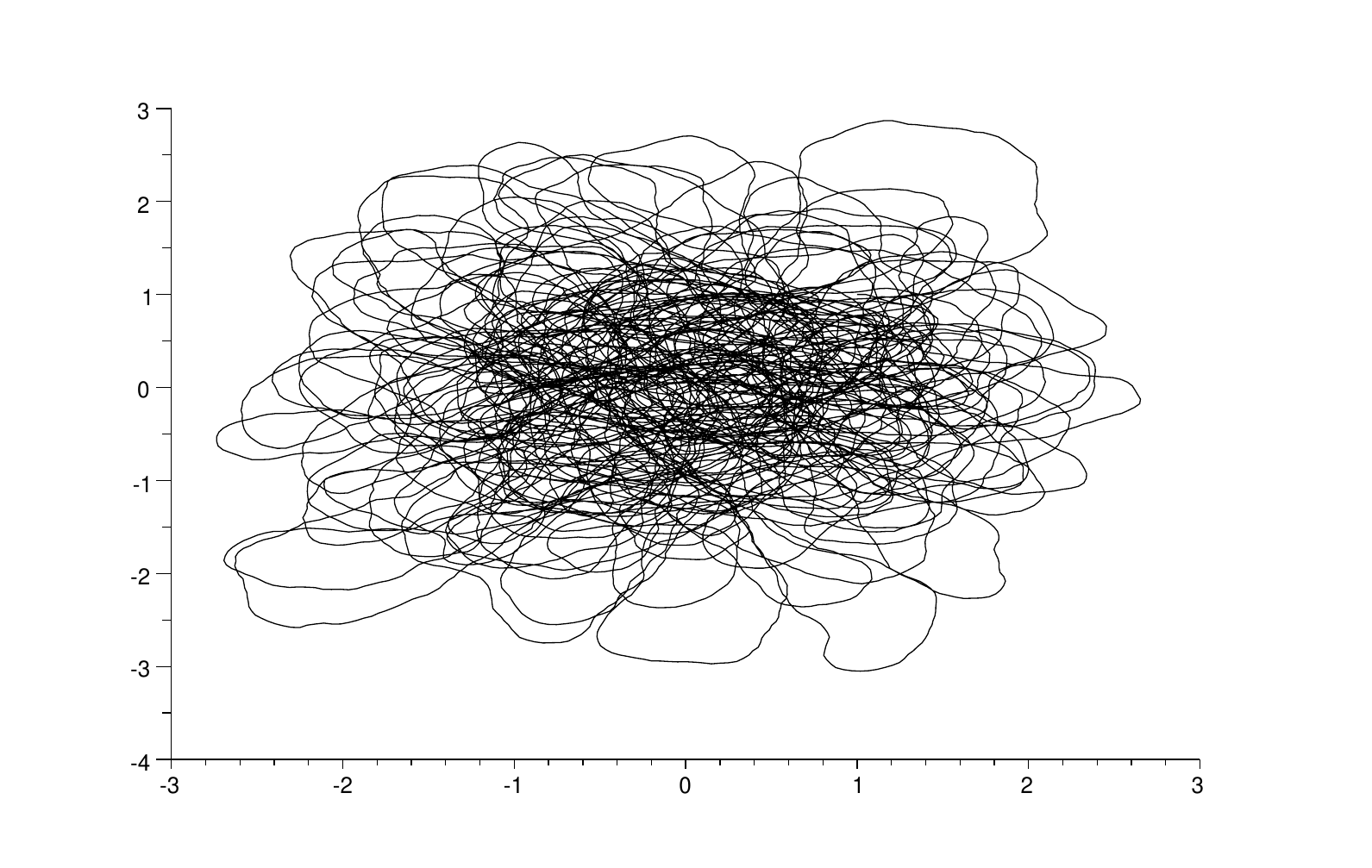}}
\captionof{figure}{\scriptsize{Evolution of the coordinate $(U_{t},V_{t})$ after time $T=750$, where $\sigma$ is respectively $0.1, 1$ and $4$} }
\end{center}
\begin{center}
\parbox{0.3\linewidth}{\includegraphics[width=0.35\textwidth]{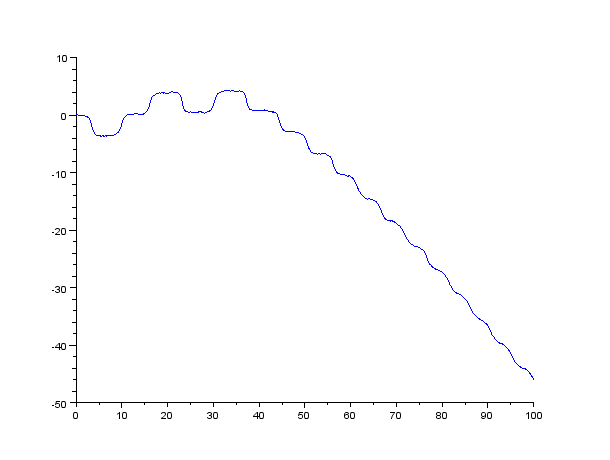}}
\parbox{0.2\linewidth}{}
\parbox{0.3\linewidth}{\includegraphics[width=0.35\textwidth]{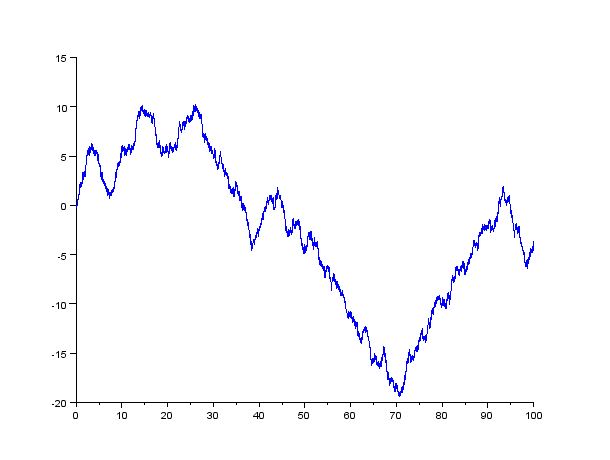}}
\parbox{0.2\linewidth}{}
\parbox{0.3\linewidth}{\includegraphics[width=0.35\textwidth]{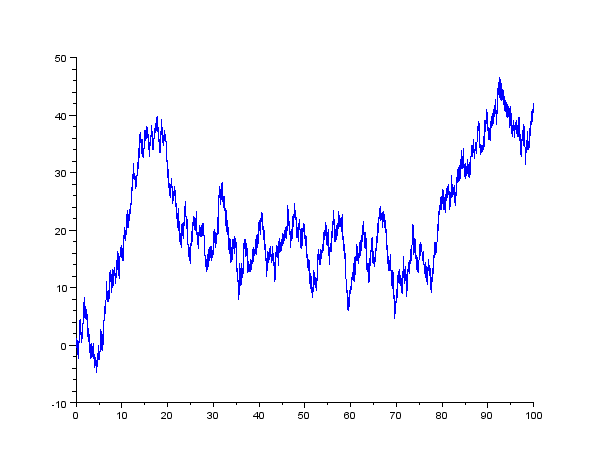}}%
\captionof{figure}{\scriptsize{Evolution of the angle $X_{t}$ after time $T=100$, where $\sigma$ is respectively $0.1, 1$ and $4$.} }
\end{center}
\subsubsection*{\textbf{The zero noise limit}}
We point out that (\ref{eq00}) is -for $\sigma\ll 1$- a random perturbation of the following ordinary differential equation (ODE)
\begin{equation}\label{eq001}
\left\{
\begin{array}{r c l}
\dot{X}_{t}&=& \sin(X_{t})U_{t}-\cos(X_{t})V_{t}\\
\dot{U}_{t}&=& \cos(X_{t}) \\
\dot{V}_{t}&=& \sin(X_{t})
\end{array}
\right.
\end{equation}
\begin{center}
\parbox{0.4\linewidth}{\includegraphics[width=0.40\textwidth]{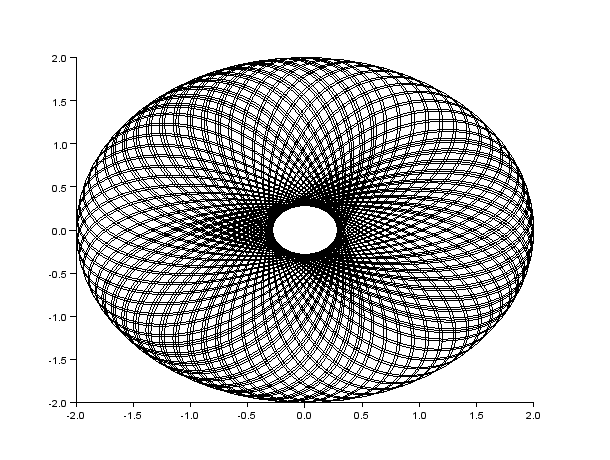}}
\parbox{0.2\linewidth}{}
\parbox{0.4\linewidth}{\includegraphics[width=0.40\textwidth]{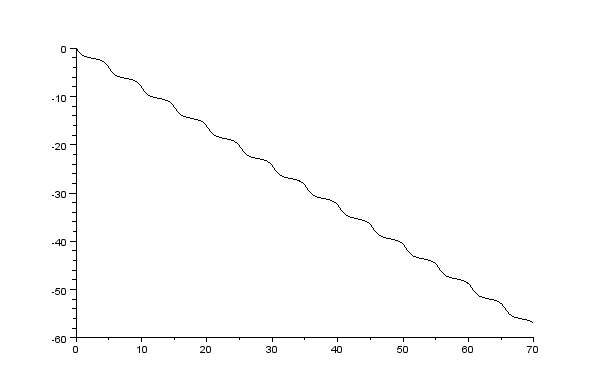}}
\captionof{figure}{\scriptsize{Evolution of $(U_{t},V_{t})$ after time $T=1000$ (left) and evolution of $X_{t}$ until time $T=70$ (right). Both simulations started with initial condition $(x,u,v)=(0,0,2)$.}}
\end{center}

The dynamics of (\ref{eq001}) can be fully described as follows:

Let ${\cal R} : \mathbb{S}^1 \times \mathbb{R}^2 \mapsto \mathbb{S}^1 \times \mathbb{R}^2$ be the map defined by
${\cal R}(x,z) = (x, {\cal R}_x z)$ where ${\cal R}_x$  is the rotation of angle $x.$  Let $H :   \mathbb{R}^2 \mapsto [1, \infty] $ be the map defined by
$$H(u,v) = \left\{ \begin{array}{l}
             \frac{1}{2}(u^2 + v^2 - \log(v^2)),  \mbox{ if } v \neq 0,\\
             \infty, \mbox{ if } v = 0.
           \end{array}  \right.$$
Set $H_c = H^{-1}(c).$  Then  $H_{\infty}$ is the line $v = 0,$ while for $c < \infty,$  $H_c$ has two components $H^+_c$ and $H^{-}_c$ obtained from each other by reflection along the line $v = 0.$
For $c > 1/2$, $H_c^+$ is  a closed curve around $(0,1),$  and $H_{1/2}^+ = \{(0,1)\}.$

Given $\alpha \in \{-,+\}$ and $c \in [1,\infty[$ set $\mathbb{T}^{\alpha}_c = {\cal R} (\mathbb{S}^1 \times H^{\alpha}_c)$ and $\mathbb{T}_{\infty} = {\cal R} (\mathbb{S}^1 \times H_{\infty}).$ Then  $\mathbb{T}^{\alpha}_{1/2}$ is a closed curve, $\mathbb{T}_c^{\alpha}$ is, for $c > 1$, a torus and $\mathbb{T}_{\infty}$ is a full twisted strip. Furthermore
\begin{equation}
\label{foliation}
\mathbb{S}^1 \times \mathbb{R}^2 = \bigcup_{c \geq 1, \alpha \in \{-,+\}} \mathbb{T}^{\alpha}_c \cup \mathbb{T}_{\infty}
\end{equation}
\begin{theorem}\label{theoremCercle}
The foliation (\ref{foliation}) is invariant under the dynamics (\ref{eq001}). More precisely,
\begin{description}
\item (i)
$\mathbb{T}^{\alpha}_{1/2}$ consists  of a periodic orbit having  period $2\pi;$
\item (ii) For $c > 1/2$  the orbits on $\mathbb{T}_c^{\alpha}$ are either all periodic or all dense in  $\mathbb{T}^{\alpha}_c$. Furthermore, the set of $c$ such that the orbits on $\mathbb{T}_c^{\alpha}$ are periodic is a countable and dense subset of $]1/2,\infty [$ ;
\item (iii) On $\mathbb{T}_{\infty}$ the solution to (\ref{eq001}) with initial condition $(x_0,u_0,v_0)$ is given by $$(X_{t},U_{t},V_{t}) = (x_0, u_0 + t \cos(x_0), v_0 + t \sin(x_0)).$$
\end{description}
\end{theorem}
The proof is the purpose of the appendix.
\begin{rem} To determine whether or not the orbits on $\mathbb{T}_c^{\alpha}$ are periodic, we introduce (see appendix) some function $T_{.,2}: ]1/2,\infty [\rightarrow ]2\sqrt{2},\sqrt{2}\pi[$ which is continuous and decreasing and prove that the orbits on $\mathbb{T}^{\alpha}_c$ are periodic if and only if $\frac{T_{c,2}}{2\pi}\in\mathbb{Q}$. Details are given in the appendix.
\end{rem}
\begin{flushleft}
\parbox{0.4\linewidth}{\includegraphics[width=0.60\textwidth]{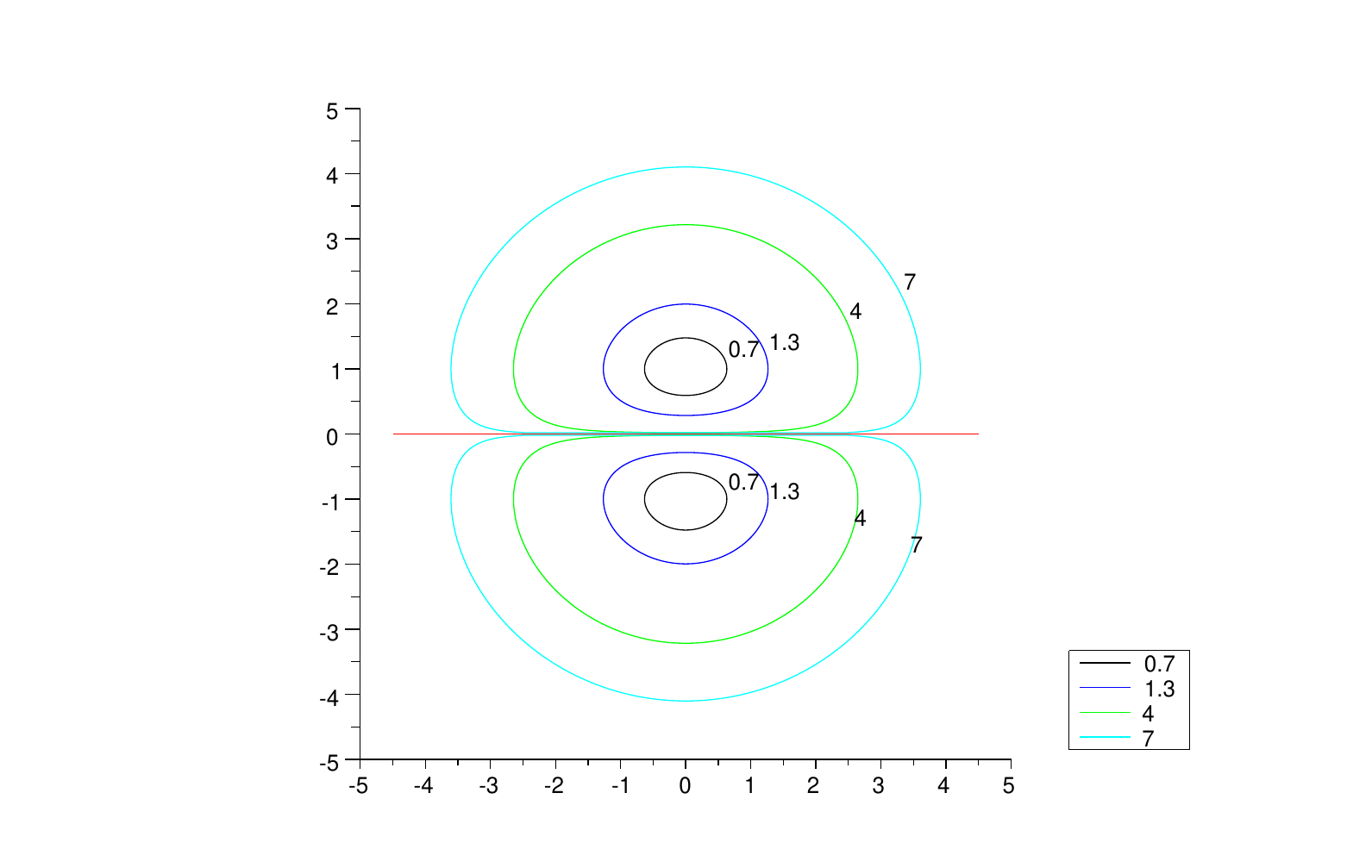}}
\parbox{0.4\linewidth}{\includegraphics[width=0.50\textwidth]{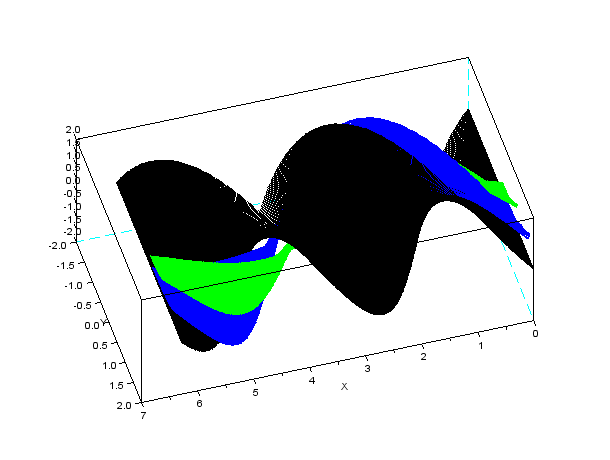}}
\captionof{figure}{\scriptsize{The left picture shows level sets of the function $H$ whereas the right picture shows the full twisted strip (in black) and two torus $T_{c}^{+}$, with $c=\sqrt{2}$ (in green) and $c= 2$ (in blue).}}
\end{flushleft}
\section{Description of the model}\label{s2}
Let us start by fixing some notation. Throughout all the paper, we let $\nabla$ denote the gradient on $M$, $\Delta_{M}$ the Laplacian on $M$ and for some vector field $\mathcal{X}$ on a manifold $\mathcal{N}$, we denote by $\mathcal{X}(f)$ the Lie derivative of $f$ along $\mathcal{X}$; $f$ being a smooth function.
\medskip

For a smooth function $V:M\times M\rightarrow \mathbb{R}$ and for a Borel measure $\mu$, we let $V\mu:M\rightarrow\mathbb{R}$ denotes the function defined by $$V\mu (x)=\int_{M}V(u,x)\mu(du). $$
We then consider the model
\begin{equation}\label{eq1}
	dX_{t}= \sigma\sum_{j=1}^{N}F_{j}(X_{t})\circ dB_{t}^{(j)}-\nabla V\mu_{t}(X_{t})dt,\qquad X_{0}=x,
\end{equation}

 where $\sigma >0$, $(B^{(1)},\cdots,B^{(N)})$ is a standard Brownian motion on $\mathbb{R}^{N}$, $\circ$ denotes the Stratonovitch integral, $\{ F_{i}\}$ is a family of smooth vectors fields on $M$ such that $$ \sum_{i=1}^{N}F_{i}(F_{i}f)=\Delta_{M}f,\: f\in\mathcal{C}^{\infty}$$
and $\mu_{t}$ is the random occupation measure defined by
$$\mu_{t}=\int_{0}^{t}\delta_{X_{s}}ds .$$
Note that there exists at least one such family $\{F_{i}\}$ since by Nash's embedding Theorem, there exists $N\in\mathbb{N}$ large enough such that $M$ is isometrically embedded in $\mathbb{R}^{N}$ with the standard metric (see Theorem 3.1.4 in \cite{HSU} or Proposition 2.5 in \cite{BRL}).\\
In this paper, we suppose that the function $V$ has the following form
\begin{equation}\label{eq2} V(x,y)=\sum_{j=1}^{n}a_{j}e_{j}(x)e_{j}(y),
\end{equation}
where $(e_{j})_{j=1,\cdots,n}$ are eigenfunctions for the Laplacian associated to non zero eigenvalues $\lambda_{1},\cdots,\lambda_{n}<0$ such that $$\int_{M}e_{j}(x)e_{k}(x)dx=\delta_{k,j},$$ where $\delta_{k,j} $ is the Kronecker symbol and $dx$ stands for the Riemannian measure on $M$. We also assume that $a_{j}>0$ for all $j=1,\cdots,n$.\\
Due to the particular form for $V$, we can obtain a "true" stochastic differential equation by introducing the new variables $U_{k,t}=\int_{0}^{t}e_{k}(X_{s})ds$. Therefore we get the following system on $\mathbb{M}:=M\times\mathbb{R}^{n}$
\begin{equation}\label{eq3}
\left\{
\begin{array}{r c l}
dX_{t}&=&\sigma\sum_{j=1}^{N}F_{j}(X_{t})\circ dB_{t}^{(j)}-\sum_{j=1}^{n}a_{j}\nabla e_{j}(X_{t})U_{j,t}dt\\
dU_{k,t}&=&e_{k}(X_{t})dt, \qquad\qquad k=1,\cdots,n
\end{array}
\right.
\end{equation}
with initial condition $(x,0,\cdots,0)$. 
In the rest of the paper, we will work with the system (\ref{eq3}) and prove that:
\begin{enumerate}
\item There exists a unique global strong solution for the system (\ref{eq3});
\item Strong Feller property holds;
\item The system admits a unique invariant measure which is given explicitly as the product of the uniform probability on $M$ and a Gaussian probability on $\mathbb{R}^{n}$;
\item The law of the solution converges to $\mu$ exponentially fast.
\end{enumerate}
The paper is organized as follows. In the next section, we present the main results and the proof of point 1.

In Section 4, we provide the proofs of points 2 and 3. To this end, we introduce a property, called \textit{condition $(E')$} and prove that it implies the Strong Feller property.

In section 5 is given the proof of an exponential decay in $L^{2}(\mathbb{M},\mu)$, where $\mu$ is the unique invariant probability whereas a proof for an exponential decay in the Total Variation norm is presented in Section 6.
\section{Presentation of the results}\label{result}
Recall that $\mathbb{M}=M\times\mathbb{R}^{n}$. Throughout, we denote by $\mathcal{C}_{0}(\mathbb{M})$ the set of function $f:\mathbb{M}\rightarrow \mathbb{R}:(x,u)\mapsto f(x,u)$ which are continuous and such that  $f(x,u)\rightarrow 0 $ when  $\Vert u\Vert\rightarrow\infty$,  and by $\mathcal{C}_{c}^{k}(\mathbb{M})$ the set of function which are $k$ times continuously differentiable with compact support.\\
We equip $\mathcal{C}_{0}(\mathbb{M})$ with the supremum norm
$$\Vert f\Vert_{\infty}:=sup_{y\in \mathbb{M}}\vert f(y)\vert. $$

Let $G_{0},\: G_{1},\cdots,\: G_{N}$ be the vector fields on $\mathbb{M}$ defined by
$$G_{0}(x,u)=\begin{bmatrix}
   -\sum_{j=1}^{n}a_{j}\nabla e_{j}(x)u_{j}\\
   e_{1}(x) \\
	\vdots\\
   e_{n}(x)
\end{bmatrix},
 $$
and for $j=1,\cdots,N$,
$$G_{j}(x,u)=\begin{bmatrix}
   \sigma F_{j}(x)\\
   0 \\
	\vdots\\
   0
\end{bmatrix},
 $$
with $x\in M$ and $u\in\mathbb{R}^{n}$. So (\ref{eq3}) can be rewritten as:
\begin{equation} \label{eq5}
dY_{t}=\sum_{j=1}^{N}G_{j}(Y_{t})\circ dB_{t}^{j} +G_{0}(Y_{t})dt.
\end{equation}
\begin{prop}\label{p0}
For all $y=(x,u)\in\mathbb{M}$ there exists a unique global strong solution $(Y_{t}^{y})_{t\geqslant 0}$ to (\ref{eq5}) with initial condition $Y_{0}^{y}=y=(x,u)$. Moreover, we have
\begin{equation}\label{eq7}
 Y_{t}^{y}=(X_{t}^{y},U_{t}^{y})\in M\times \bar{B}(u,K t),
\end{equation}
where $K=(\max_{y\in M} \sum_{j=1}^{n}e_{j}(y)^{2})^{1/2}$ and $\bar{B}(u,R)=\{ v\in\mathbb{R}^{n}\; :\; \Vert v-u\Vert\leqslant R\}$.
\end{prop}
\begin{proof1} Existence and uniqueness is standard since $G_{0}$ is locally Lipschitz and sub-linear (see for example \cite{RY}, page 383). 
Concerning (\ref{eq7}), note that we have
\begin{eqnarray*}
\sum_{j=1}^{n}(U_{j,t}-u_{j})^{2}&\leqslant & t\int_{0}^{t}\sum_{j=1}^{n}e_{j}(X_{s})^{2}ds\\
&\leqslant & t^{2} \max_{y\in M} \sum_{j=1}^{n}e_{j}(y)^{2}<\infty;
\end{eqnarray*}
which proves (\ref{eq7}).
\end{proof1}
Throughout, we let $(P_{t})_{t\geqslant 0}$ denote the semi-group induced by (\ref{eq5}). Recall that for any bounded or nonnegative measurable function $f:\mathbb{M}\rightarrow \mathbb{R}$, $P_{t}f$ is the function defined by 
\begin{equation}\label{SemGroup} P_{t}f(y)=\mathbb{E}(f(Y_{t}^{y}))\text{ for all } y\in\mathbb{M}.\end{equation} 
\begin{lem} \label{l9} The semi-group $(P_{t})_{t\geqslant 0}$ is Feller, meaning that
\begin{enumerate}
\item For all $t\geqslant 0$, $P_{t}(\mathcal{C}_{0}(\mathbb{M}))\subset \mathcal{C}_{0}(\mathbb{M})$.
\item For all $f\in\mathcal{C}_{0}(\mathbb{M})$, $\lim_{t\rightarrow 0}\Vert P_{t}f -f\Vert_{\infty}=0$.
\end{enumerate}
\end{lem}
\begin{proof1} By Proposition \ref{p0}, for all $T>0$, $(Y_t^y)_{t\in [0,T]}$ lies on a deterministic compact set depending only $y$ and $T$. Hence, by standard results (see eg Theorem IX.2.4 in \cite{RY})
, $y\mapsto Y_{t}^{y}$ is continuous. Thus, by dominated convergence, $y\mapsto P_{t}f(y)$ lies in $\mathcal{C}_{0}(\mathbb{M})$ for all $f\in \mathcal{C}_{0}(\mathbb{M})$.\\
In order to prove the second point, it suffices to show that $\lim_{t\downarrow 0}P_{t}f(y)=f(y)$ (see Proposition III.2.4 in \cite{RY}). This follows again from continuity of $t\mapsto Y_{t}^{y}$ and dominated convergence.
\end{proof1}
The next result gives further informations on the semi-group.
\begin{prop}\label{p5} The set $\mathcal{C}_{c}^{2}(\mathbb{M})$ is stable for $P_{t},\: t\geqslant 0$, ie for all $t\geqslant 0$, $P_{t}(\mathcal{C}_{c}^{2}(\mathbb{M}))\subset \mathcal{C}_{c}^{2}(\mathbb{M})$.
\end{prop}
\begin{proof1} Let $f\in\mathcal{C}_{c}^{2}(\mathbb{M})$. The fact that $P_{t}f$ has a compact support is a consequence of Equation (\ref{eq7}). Let us now prove that $P_{t}f$ is twice continuously differentiable.\\
Let $y=(x_{0},u)\in\mathbb{M}$ and $R>0$. For $\tilde{y}\in M\times B(u,R)$, we have, by Proposition \ref{p0},
\begin{equation}
(Y_{s}^{\tilde{y}})_{0\leqslant s\leqslant t}\subset M\times \bar{B}(u,Kt+R).
\end{equation}
Pick a smooth function $\psi:\mathbb{R}^{n}\rightarrow\mathbb{R}_{+}$ which is 1 on the ball $B(u,Kt+R)$, 0 outside the ball $\bar{B}(u,Kt+R+1)$ and $\psi (v)\leqslant 1$ for all $v$.\\
Consider now the SDE defined by
\begin{equation}\label{EqYtilde}
d\tilde{Y}_{t}=\sum_{j=1}^{N}G_{j}(\tilde{Y}_{t})\circ dB_{t}^{j} +\tilde{G}_{0}(\tilde{Y}_{t})dt,
\end{equation}
where $\tilde{G}_{0}(x,v)=G_{0}(x,u+\psi (v)(v-u))$. Let us denote by $\tilde{P}_{t}$ its associated semi-group. The fact that $G_{0}$ is smooth and locally Lipschitz implies that $\tilde{G}_{0}$ is smooth and Lipschitz. By Nash's embedding Theorem and proceeding in the same way as in Proposition 2.5 in \cite{BRL}, we can extend (\ref{EqYtilde}) to a SDE on $\mathbb{R}^N\times \mathbb{R}^n$ and $f$ to a function in $\mathcal{C}^2(\mathbb{R}^N\times \mathbb{R}^n)$.
Therefore, in view of subsection 3.2.1 in \cite{DaPrKOLMOG} and of Proposition 2.5 in \cite{DaP-Rock}, it follows that $\tilde{P}_{s}f $ is a function of class $\mathcal{C}^{2}$ for all $s\geqslant 0$. Since
\begin{equation}
P_{s}f(\tilde{y})=\tilde{P}_{s}f(\tilde{y})\: \text{ for all } 0\leqslant s\leqslant t\: \text{ and all } \tilde{y}\in M\times B(u,R),
\end{equation}
it follows that $P_{t}f$ is of class $\mathcal{C}^{2}$ on $M\times B(u,R)$.\\
Consequently, $P_{t}f\in\mathcal{C}_{c}^{2}(\mathbb{M})$.
\end{proof1}
The infinitesimal generator of $(P_{t})_{t\geqslant 0}$ is the operator
\begin{equation}\label{GenINFIN} \mathcal{L}: D(\mathcal{L})\rightarrow \mathcal{C}_{0}(\mathbb{M}): f\mapsto \lim_{t\downarrow 0} \frac{P_{t}f-f}{t},
\end{equation}
where $D(\mathcal{L}):=\{ f\in \mathcal{C}_{0}(\mathbb{M}):\; \frac{P_{t}f-f}{t} \text{ converges in }\mathcal{C}_{0}(E) \text{ when } t\downarrow 0 \}$.
Then (see for example Theorem 17.6 in \cite{Kallenberg}) for all $f\in D(\mathcal{L})$,
\begin{equation}\label{eq8}
P_{t}f -f=\int_{0}^{t}\mathcal{L}(P_{s}f)ds=\int_{0}^{t}P_{s}(\mathcal{L}f)ds
\end{equation}
We briefly recall the following result which characterize the elements of $D(\mathcal{L})$:
\begin{theorem} (Propositions VII.1.6 and VII.1.7 in \cite{RY})\label{t1}\\
For $g,h\in\mathcal{C}_{0}(\mathbb{M})$, the following assertions are equivalent:
\begin{enumerate}
\item $h\in D(\mathcal{L})$ and $\mathcal{L}h=g$.
\item For all $y\in E$, the process $$h(Y_{t}^{y})-\int_{0}^{t}g(Y_{s}^{y})ds $$ is a martingale with respect to the filtration $\mathcal{F}_{t}=\sigma (Y_{s}^{y}\; :\; 0\leqslant s\leqslant t)$.
\end{enumerate}
\end{theorem}
Since the definition of the infinitesimal generator is implicit, it is convenient to introduce a more tractable operator: the \textit{Kolmogorov operator}.
\begin{defin} The Kolmogorov operator associated to (\ref{eq3}) is the operator defined on $\mathcal{C}^{2}$ bounded functions having first and second bounded derivatives by
$$L= \frac{\sigma^{2}}{2}\Delta_{M}-\sum_{k=1}^{n}a_{k}u_{k}(\nabla e_{k}(x),\nabla_{x} \stackrel{ \textbf{.}}{})_{TM}+\sum_{k=1}^{n}e_{k}(x)\partial_{u_{k}}, $$
with the convention $(\Delta_{M}f)(x,u)=(\Delta_{M}f(.,u))(x)$ and $(.,.)_{TM}$ stands for the inner product on the tangent bundle of $M$.\\
\end{defin} 
The link between the infinitesimal and the Kolmogorov operator is given by the next proposition.
\begin{prop}\label{p5b} Let $f$ be a $\mathcal{C}^{2}$ bounded function having first and second bounded derivatives, then $f\in D(\mathcal{L})$ and
$$Lf=\mathcal{L}f.$$
\end{prop}
\begin{proof1} It follows from It\^o's formula and Theorem \ref{t1}.
\end{proof1}

\begin{defin}\label{probinv} Let $\Phi:\mathbb{R}^{n}\rightarrow \mathbb{R}:u=(u_{1}\cdots,u_{n})\mapsto \ln(C(\Phi))+\frac{1}{2}\sum_{k=1}^{n}a_{k}\vert\lambda_{k}\vert u_{k}^{2}$ with $$C(\Phi)=\int_{\mathbb{R}^{n}}\exp(-\frac{1}{2}\sum_{k=1}^{n}a_{k}\vert\lambda_{k}\vert u_{k}^{2})du=\prod_{i=1}^{n}\sqrt{\frac{2\pi}{\vert\lambda_{i}\vert a_{i}}}<\infty.$$ Recall that $\lambda_{i}<0$ is the eigenvalue associated to the eigenfunction $e_{i}$ of $\Delta_{M}$.
 On $\mathbb{M}$, we define the probability measure \begin{equation}\mu(dx\otimes du)=\nu(dx)\otimes e^{-\Phi(u)}du=:\varphi(y)dy, \end{equation}
with $y=(x,u)$ and $\nu(dx)=\frac{dx}{\int_{M}dz}$ is the uniform probability measure on $M$.
\end{defin}
\begin{rem} \label{r2} Note that $\mu(dy)$ does not depend on the noise term $\sigma$.
\end{rem}
We can now state our first main result.
\begin{theorem} \label{t2} Let $(P_{t})_{t\geqslant 0}$ be the semi-group associated to the system (\ref{eq3}) and $P_{t}(y_{0},dy)$ its transition probability. Then
\begin{enumerate}
\item[1)] The semi-group $(P_{t})_{t\geqslant 0}$ is strongly Feller (meaning that $P_{t}f$ is a bounded continuous function for whatever bounded measurable function $f$) and there exists a $\mathcal{C}^{\infty}((0,\infty),\mathbb{M},\mathbb{M})$ function $p_{t}(y_{0},y)$ such that $P_{t}(y_{0},dy)=p_{t}(y_{0},y)dy$ for all $y_{0}\in\mathbb{M}$ and $(L_{z}^{\ast}-\partial_{t})p_{t}(y,z)=0$,
\item[2)] The probability $\mu (dy)=\varphi(y)dy$, where $\varphi$ is given in Definition \ref{probinv}, is the unique invariant probability. 
 Moreover for all $y\in\mathbb{M}$ and for all bounded measurable function $f$, we have $$\lim_{t\rightarrow\infty}P_{t}f(y)=\int_{\mathbb{M}}f(z)\mu(dz). $$
 Furthermore, the process $(Y_{t})_{t}$ is positive Harris recurrent, ie for all Borelian set $R$ such that $\mu (R)>0$, then $$\int_{0}^{\infty}\mathbf{1}_{R}(Y_{t}^{y})dt=\infty\: a.s$$ for all $y\in\mathbb{M}$.
\item[3)]  $\lim_{t\rightarrow\infty}\int_{\mathbb{M}}\vert p_{t}(z,y)-\varphi(y)\vert dy =0$ for all $z\in \mathbb{M}$.
\end{enumerate}
\end{theorem}
\begin{rem} The fact that $\mu$ is independent of the parameter $\sigma$ implies that it is also an invariant probability of the deterministic system obtained with $\sigma =0$. However, in that case it is not necessarily unique (compare with Theorem \ref{theoremCercle}, where there exists infinitely many compact disjoint invariant sets, thus infinitely many ergodic probabilities.)
\end{rem}

As an immediate consequence of the Harris positive recurrence property, we have
\begin{cor} \label{corErg} For all $f\in L^{1}(\mu)$, $$\frac{1}{t}\int_{0}^{t}f(Y_{s}^{y})ds\rightarrow \int_{\mathbb{M}}f(y)\mu(dy)$$  almost surely for any $y\in\mathbb{M}$.
\end{cor}
\begin{proof1} Apply Theorem 3.1 in \cite{A-D-R} to the positive and negative part of $f$. 
\end{proof1}
The next results establish exponential rate of convergence of $(P_{t})_{t\geqslant 0}$ to $\mu$.
\begin{theorem} \label{t4} For every $\eta >0$ and $g\in L^{2}(\mu)$
$$\Vert P_{t}g-\int_{\mathbb{M}}g(y)\mu (dy)\Vert_{L^{2}(\mu)}\leqslant \sqrt{1+ 2\eta}\Vert g-\int_{\mathbb{M}}g(y)\mu (dy)\Vert_{L^{2}(\mu)}e^{-\lambda t}, $$
where $$\lambda=\frac{\eta}{1+\eta}\frac{K_{1}\sigma^{2}}{1+K_{2}\sigma^{2}+K_{3}\sigma^{4}}\:,$$
with $$K_1 =\frac{1}{4(2+(1+N_2)^2)}(\frac{\Lambda}{1+\Lambda})^2 , $$
$$K_2 =\frac{(1+N_2)\sum_{j=1}^n\vert \lambda_j\vert}{2+(1+N_2)^2}, $$ 
$$K_3 = \frac{(\sum_{j=1}^n\vert \lambda_j\vert)^2}{4(2+(1+N_2)^2)}, $$
$$\Lambda =\min_{i=1,\ldots , n}\vert \lambda_i\vert a_i $$ and 
$$N_2=2 \frac{n}{\min\{\vert \lambda_j\vert , j=1,\cdots,n\}}\sup_{i = 1, \ldots, n} \|\nabla e_i\|_{\infty}^2  \sqrt{4+\sum_{i=1}^n |\lambda_i| a_i} + 4 \|\sum_i e_i^2\|_{\infty}. $$ 
\end{theorem}
\begin{rem}\label{rem5} Note that if $g\in L^{2}(\mu)$, then it is not clear at first glance that $P_{t}g$ is meaningful. However it is. In order to prove it, set $h_{t}(y,z)=p_{t}(y,z)/\varphi (z)$. Due to the properties of $p_{t}(y,.)$ and $\varphi$ for all $t>0$ and $x\in\mathbb{M}$ (see Theorem \ref{t2}, Proposition \ref{p0} and Definition \ref{probinv}), then $h_{t}(y,.)$ has compact support. Thus, by the Cauchy-Schwarz inequality, we obtain 
\begin{equation}\mathbb{E}(\vert g\vert (Y_{t}^{y}))=\int_{\mathbb{M}}\vert g\vert (z)p_{t}(y,z)dz=\int_{\mathbb{M}}\vert g\vert (z)h_{t}(y,z)\mu(dz)\leqslant\Vert g \Vert_{L^{2}(\mu)}\Vert h_{t}(y,.)\Vert_{L^{2}(\mu)}.\end{equation}
Furthermore, we have $P_{t}g\in L^{2}(\mu)$. Indeed by Jensen inequality and invariance of $\mu$, we have $\int_{\mathbb{M}}(P_{t}g)^{2}(y)\mu(dy)\leqslant \int_{\mathbb{M}}P_{t}(g^{2})(y)\mu(dy)=\int_{\mathbb{M}}g^{2}(y)\mu(dy)<\infty$.
\end{rem}
Since both $\mu(dy)$ and $P_{t}(y_{0},dy)$ have smooth densities with respect to the Lebesgue measure for all $y_{0}\in\mathbb{M}$ and in view of the third point of Theorem \ref{t2}, we would hope to get a convergence speed for the total variation norm. Once again the answer is positive as shown by the following theorem.
\begin{theorem}\label{t5}For all $z_{0}\in\mathbb{M}$ and $t\geqslant 1$,
 $$\Vert P_{t}(z_{0},dz)-\mu(dz)\Vert_{TV}\leqslant \sqrt{1+ 2\eta}\Vert h(1,z_{0},z)-1\Vert_{L^{2}(\mu)}e^{-\lambda (t-1)}, $$
where $ h(1,z_{0},z)=\frac{p_{1}(z_{0},z)}{\varphi(z)}$ and $$\lambda=\frac{\eta}{1+\eta}\frac{K_{1}\sigma^{2}}{1+K_{2}\sigma^{2}+K_{3}\sigma^{4}}\:,$$
$\mu$ is the probability given in Theorem \ref{t2} and the constants $K_{j}<\infty,\; j=1,2,3$, are the same as in Theorem \ref{t4}.
\end{theorem}
The proofs of Theorem \ref{t4} and \ref{t5} are postponed to Sections \ref{SectL2cvge} and \ref{SectTVcvge}. 
\section{Proof of Theorem \ref{t2}}
We emphasize, from Equation (\ref{eq5}), that the Kolmogorov operator $L$ can be expressed in H\"ormander's form (as a sum of squares):
\begin{equation}L=\frac{1}{2}\sum_{j=1}^{N}G_{j}^{2}+G_{0}, \end{equation}
where $G_{j}^{2}(f)=G_{j}(G_{j}f)$.
The proof mainly relies on classical results by Kanji Ichihara and Hiroshi Kunita in \cite{Ichihara} dealing with this type of operator.
\subsection*{\textbf{Proof of assertion 1): the Strong Feller Property.}}
Throughout, we use the following notation. If  ${\cal N}$ is a smooth manifold (such as $M, \mathbb{M}$ or $\RR^m$), $W : {\cal C}^{\infty}({\cal N}) \rightarrow {\cal C}^{\infty}({\cal N})$ a linear map (typically a differential operator) and  $f : {\cal N} \rightarrow \RR^n : x \mapsto (f_1(x), \ldots, f_n(x))$  a smooth map, we let $W(f) : {\cal N} \rightarrow \RR^n$ denote the map defined by $$W(f)(x) = (W(f_1)(x), \ldots, W(f_n)(x)).$$
Given two smooth vector fields $A$ and $B$ on  ${\cal N}$ recall that the {\em Lie-bracket} of $A$ and $B$ is the vector field on ${\cal N}$ characterized by
$$[A,B](f)=A(B(f))-B(A(f))$$ for all $f \in \mathcal{C}^{\infty}({\cal N}).$
In case ${\cal N} = \RR^m$ then for all $x \in \RR^m$
$$[A,B](x)=DB(x)A(x)-DA(x)B(x)$$
where $DA(x)$ (resp. $DB(x)$ ) stands for the derivative of $A$ (resp. $B$) at $x.$

Let ${\cal G}_{0}=\{ G_{1},\cdots,G_{N}\}$. Define then recursively ${\cal G}_{k}, k\geqslant 1,$ by
\begin{eqnarray*}\label{DEFinductive}{\cal G}_{k}={\cal G}_{k-1}\cup \{[B,G_{j}],\; B\in {\cal G}_{k-1}\text{ and } j=0,\cdots,N.\} \end{eqnarray*}
Let then ${\cal G}_{\infty}=\bigcup_{k\geqslant 0}{\cal G}_{k}$  and for all $(x,u) \in \mathbb{M}$
$${\cal G}_{\infty}(x,u) = \{ V(x,u) \; : V \in {\cal G}_{\infty} \}.$$
Using the terminology of \cite{Ichihara}, we say that
\begin{defin} The dynamics (\ref{eq5}) satisfies the ellipticity condition $(E)$ if for all $(x,u) \in \mathbb{M}, {\cal G}_{\infty}(x,u)$ spans $T_{(x,u)} \mathbb{M} = T_x M \times \RR^n.$
\end{defin}
The next result rephrases Lemma 5.1 $(ii)$ and Theorem 3 $(i)$ and $(iii)$ of \cite{Ichihara}.
\begin{lem}
\label{Ichilem}
If (\ref{eq3}) satisfies $(E)$ then the induced semi-group $(P_t)$ is strongly Feller and there exists a $\mathcal{C}^{\infty}((0,\infty),\mathbb{M},\mathbb{M})$ function $p_{t}(y_{0},y)$ such that $P_{t}(y_{0},dy)=p_{t}(y_{0},y)dy$ for all $y_{0}\in\mathbb{M}$ and $(L_{z}^{\ast}-\partial_{t})p_{t}(y,z)=0$.
\end{lem}
\begin{rem}
Note that when $\sigma =0$, the condition $(E)$ is never satisfied since ${\cal G}_{0}$ is reduced to $\{0\};$ hence ${\cal G}_{\infty}=\{0\}.$
\end{rem}
Let ${\cal A}_{0} = \{ F_{1},\cdots,F_{N}\}$
and for all $k \geq 1$
\begin{eqnarray}
\label{DEFinductive2} {\cal A}_{k}={\cal A}_{k-1}\cup \{F_{j} B,\; B \in {\cal A}_{k-1}\text{ and } j=1,\cdots,N\},\end{eqnarray}
where $F_{j}B$ is the operator on $\mathcal{C}^{\infty}(M)$ defined by $(F_{j}B)(f)=F_{j}(B(f))$.

Let then ${\cal A}_{\infty}=\bigcup_{k\geqslant 0}{\cal A}_{k}$  and for all $x \in M$
$${\cal A}_{\infty}(x) = \{ W(e)(x) \; : W \in {\cal A}_{\infty} \}$$
where $e : M \rightarrow \RR^n$ is the map defined by $e(x) = (e_1(x), \ldots, e_n(x)).$
Note that while ${\cal G}_{\infty}$ is a set of vector fields on $\mathbb{M},$ ${\cal A}_{\infty}$ is a set of differential operators of all orders on $\mathcal{C}^{\infty}(M).$
\begin{defin} We say that the condition $(E')$ is fulfilled if and only if for all $x\in M$,
${\cal A}_{\infty}(x)$ spans $\RR^n.$
\end{defin}

\begin{lem}\label{EeE} Suppose $\sigma > 0.$ Then, condition $(E')$ implies condition $(E)$.
\end{lem}
The proof relies on the following lemma.
\begin{lem}\label{l1} Let $\mathfrak{e}:\mathbb{R}^{m}\rightarrow \mathbb{R}^{n}$ be a smooth function and let $F(x,u)=\begin{bmatrix}
   A(x)\\
   0
\end{bmatrix}$ and $G(x,u)=\begin{bmatrix}
   B(x,u)\\
   \mathfrak{e} (x)
\end{bmatrix}$ be two vector fields on $\mathbb{R}^{m+n}$, where $A:\mathbb{R}^{m}\rightarrow\mathbb{R}^{m}$ and $B:\mathbb{R}^{m+n}\rightarrow \mathbb{R}^{m}$ are smooth functions. Then
$$[F,G](x,u)=\begin{bmatrix}
   [A,B(.,u)](x)\\
   A(\mathfrak{e} )(x)
\end{bmatrix},$$
with $B(.,u):\mathbb{R}^{m}\rightarrow\mathbb{R}^{m}:x\mapsto B(x,u) $
\end{lem}
\begin{proof1}  Let $(x,u)\in\mathbb{R}^{m}\times\mathbb{R}^{n}$. We then get that
$$DF(x,u)=\begin{bmatrix}
   DA & 0\\
   0 &0
\end{bmatrix}(x,u) $$
and
$$DG(x,u)=\begin{bmatrix}
   D_{x}B & D_{u}B\\
   D\mathfrak{e} &0
\end{bmatrix}(x,u). $$
Hence \begin{eqnarray}
[F,G](x,u)&=&DG(x,u)F(x,u)-DF(x,u)G(x,u)\notag\\
&=&
\begin{bmatrix}
  D_{x}B(x,u)A(x)  - DA(x)B(x,u) \\
   D\mathfrak{e}(x)A(x)
\end{bmatrix}\notag\\
&=&\begin{bmatrix}
   [A,B(.,u)](x)\\
   A(\mathfrak{e} )(x)
\end{bmatrix}
\end{eqnarray}
as stated.
\end{proof1}

\begin{proof1}[of lemma \ref{EeE}]
Let \begin{equation}\label{prodAk}W =\prod_{j=1}^{l}F_{i_{j}},\; (i_{1},\cdots,i_{l})\in\{1,\cdots,N\}^{l}.\end{equation}

By definition of $G_{0}$, and lemma \ref{l1} (used in a local chart) it follows that
\begin{equation}\label{rapCrocLie}G_{W}(x,u):=[G_{i_{1}},[\cdots ,[G_{i_{l}}, G_{0}]\cdots ]]= \sigma^{l}
\begin{bmatrix}
   \divideontimes \\
   W(e)(x)
\end{bmatrix}
\end{equation}
Thus, by hypothesis and the definition of $G_{j}$ for $j=1,\cdots ,N$,
\begin{eqnarray*} \{G_{1}(x,u),\cdots,G_{N}(x,u)\} \cup \{G_{W}(x,u) \: : W \in {\cal A}_{\infty} \} \end{eqnarray*}
spans $T_{(x,u)} \mathbb{M}$. This set being a subset of ${\cal G}_{\infty}(x,u),$ this proves the lemma.
\end{proof1}
\begin{lem}
\label{Aron}
Suppose that $\{e_1, \ldots, e_n\}$ are eigenfunctions associated to the same nonzero eigenvalue of $\Delta_M.$ Then condition $(E')$ holds true.
\end{lem}
\begin{proof1}
Let $(U, (x_1, \ldots, x_m))$ be a local chart with $U$ an open set in $M.$ Let $D_1, \ldots, D_m$ be the vector fields defined on $U$ by $D_i(f) = \frac{\partial}{\partial x_i}f.$
Define ${\cal A}_{\infty}^D$ like ${\cal A}_{\infty}$ by replacing $F_1, \ldots, F_N$ by $D_1, \ldots, D_m,$
and set ${\cal A}_{\infty}^D(x) = \{W(e)(x) \: : W \in {\cal A}_{\infty}^D\}$ for all $x \in U.$  We claim that ${\cal A}_{\infty}^D(x)$ spans $\RR^n.$
Suppose to the contrary that  there exists some $x^* \in U$ and some vector $t \in \RR^n \setminus \{0\}$ such that ${\cal A}_{\infty}^D(x^*) \subset t^\perp.$ Let $f(x) = \sum_i t_i e_i(x).$ Then $f$ is an eigenfunction of $\Delta_M$ and
for all $W \in {\cal A}_{\infty}^D$ $$ W(f) (x^*) =  W(\sum_{i = 1}^n t_i e_i)(x^*) = \langle W(e)(x^*), t \rangle = 0.$$
In other words, $f$ vanishes to infinite order at $x^*$. But
by a result of Aronzajn (see \cite{Aron}), every nonzero eigenfunction of the Laplacian on a $C^{\infty}$ manifold with  $C^{\infty}$ metric, never vanishes to infinite order.
This proves the claim.

It remains to show that ${\cal A}_{\infty}(x)$ spans $\RR^n.$ Since $F_1(x), \ldots, F_N(x)$ span $T_x M$ for all $x,$ there exist smooth real valued maps
$\alpha_{ij}, 1 \leq i \leq m, 1 \leq  j \leq N,$ defined on  $U$ such that for all $x \in U$ and $1 \leq j \leq N$
$$D_{j}(x)=\sum_{k=1}^{N}\alpha_{j,k}(x)F_{k}(x).$$
Thus
$$D_j(e)(x) = \sum_{i=1}^{N}\alpha_{j,i}(x)F_{i}(e)(x) \in span ({\cal A}_{\infty}(x)).$$
Now,  for all $\psi ,\xi\in \mathcal{C}^{\infty}(M)$ and all $H\in\mathcal{A}_{\infty}$, we have
\begin{equation*} H(\psi\xi)(x)=\psi(x) H(\xi)(x)+ \xi(x) H(\psi)(x).
\end{equation*}
Thus,
\begin{eqnarray*} D_{i}D_{j}(e)(x)
&=&\sum_{k=1}^{N}D_{i}(\alpha_{j,k})(x)F_{k}(e)(x) +\sum_{k=1}^{N}\alpha_{j,k}(x)D_{i}F_{k}(e)(x)\\
&=&\sum_{k=1}^{N}D_{i}(\alpha_{j,k})(x)F_{k}(e)(x) +\sum_{k,l=1}^{N}\alpha_{j,k}(x)\alpha_{i,l}(x)F_{l}F_{k}(e)(x) \in  span ({\cal A}_{\infty}(x))
\end{eqnarray*}
By recursion, it comes that ${\cal A}^D_{\infty}(x) \subset  span ({\cal A}_{\infty}(x))$ and since ${\cal A}^D_{\infty}(x)$ spans $\RR^n$, so does ${\cal A}_{\infty}(x)$.
\end{proof1}

\begin{lem}\label{t6} Condition $(E')$ holds.
\end{lem}
\begin{proof1}
Let $\Lambda$ be the set of distinct eigenvalues of $\{e_{1},\cdots,e_{n}\}.$  For $\lambda \in \Lambda$ let  $\{e^{\lambda}_1, \ldots, e^{\lambda}_{n(\lambda)} \} \subset \{e_{1},\cdots,e_{n}\}$ be the set of eigenfunctions having eigenvalue $\lambda$ and let $e^{\lambda} = (e^{\lambda}_1, \ldots, e^{\lambda}_{n(\lambda)}).$

Let $x\in M$. By Lemma \ref{Aron} there exist $W_{1}^{\lambda},\cdots ,W_{n(\lambda)}^{\lambda}\in\mathcal{A}_{\infty}$ such that the matrix
\begin{equation} R_{\lambda}=(W_{i}^{\lambda}(e_{j}^{\lambda})(x))_{1\leqslant i,j\leqslant n(\lambda)}
\end{equation}
has rank $n(\lambda)$.

Given a polynomial $P(x)=\sum_{j=0}^{k}\alpha_{j}x^{j}$, we let
\begin{equation}
P(\Delta_{M})=\sum_{j=0}^{k}\alpha_{j}\Delta_{M}^{j},
\end{equation}
where $\Delta_{M}^{j}$ is the operator defined recursively by $\Delta_{M}^{0}f=f$ and $\Delta_{M}^{j+1}f=\Delta_{M}^{j}(\Delta_{M}f)$ with $f\in\mathcal{C}^{2}(M)$. Note that for all $1 \leq i \leq n(\lambda)$
\begin{equation}
P(\Delta_{M})(e^{\lambda}_i)=P(\lambda)e^{\lambda}_{i}.
\end{equation}
Now let $P^{\lambda}(x)=\prod_{\substack{\alpha \in \Lambda ;}\\{\alpha \neq \lambda}}(x-\alpha)$. For $\lambda\in \Lambda$ and $i=1,\cdots ,n(\lambda)$, set
\begin{equation} H_{i}^{\lambda}= W_{i}^{\lambda}P^{\lambda}(\Delta_{M}).
\end{equation}
Then one has that $H_{i}^{\lambda}(e_{j}^{\alpha})(x)= 0$ for $\alpha\neq\lambda$ and $H_{i}^{\lambda}(e_{j}^{\lambda})(x)= P^{\lambda}(\lambda)W_{i}^{\lambda}(e_{j}^{\lambda})(x)$. Thus, the matrix $$H=(H_{i}^{\lambda}(e_{j}^{\alpha})(x))_{\lambda\in\Lambda,\; i=1,\cdots , n(\lambda)}$$ can, after a reordering if necessary, be written as a diagonal block matrix $(P^{\lambda}(\lambda)R_{\lambda}(x))_{\lambda\in\Lambda}$.

It is then easy to see that $H$ has rank $n$.
\end{proof1}
This later lemma combined with Lemmas \ref{Ichilem} and \ref{EeE} proves assertion $1)$.

\subsection*{\textbf{Proof of assertions 2) and 3). Invariant probability measure and Harris Recurrence}}
 Recall that a probability measure $\mu$ is invariant for the semi-group $(P_{t})_{t\geqslant 0}$ if
$$\int_{\mathbb{M}}P_{t}f(y)\mu(dy)= \int_{\mathbb{M}}f(y)\mu(dy)$$ for all $f\in\mathcal{C}_{0}(\mathbb{M})$.

\underline{Existence of an invariant probability measure.}
We will switch between the two notations $y\in\mathbb{M}$ and $(x,u)\in M\times\mathbb{R}^{n}$ which represent the same point. Setting
\begin{equation}\label{adjointLEB}
L^{\ast}=\frac{\sigma^{2}}{2}\Delta_{M} +\sum_{k=1}^{n}a_{k}u_{k}div_{x}(\nabla e_{k}(x) .)- \sum_{k=1}^{n}e_{k}(x)\partial_{u_{k}}.
\end{equation}
 we then observe that 
\begin{eqnarray*}
L^{\ast}\varphi(y)&=&\sum_{k=1}^{n}a_{k}u_{k}div_{x}(\nabla e_{k}(x)\varphi(y))-\sum_{k=1}^{n}e_{k}(x)\partial_{u_{k}}\varphi(y)\\
&=&\sum_{k=1}^{n}a_{k}u_{k}\lambda_{k}e_{k}(x)\varphi(y)+\sum_{k=1}^{n}e_{k}(x)a_{k}\vert\lambda_{k}\vert u_{k}\varphi(y)\\
&=&0.
\end{eqnarray*}
By Propositions \ref{p5} and \ref{p5b} together with Theorem \ref{t1}, we get for $f\in\mathcal{C}_{c}^{2}(\mathbb{M})$
\begin{eqnarray*}
\int_{\mathbb{M}}(P_{t}f(y)-f(y))\mu(dy)&=&\int_{0}^{t}\int_{\mathbb{M}}\mathcal{L}P_{s}f(y)\mu(dy)ds\\
&=& \int_{0}^{t}\int_{\mathbb{M}}LP_{s}f(y)\varphi (dy) ds
\end{eqnarray*}
Noting that for all $g,h\in\mathcal{C}_{c}^{2}(\mathbb{M})$
$$\int_{\mathbb{M}}Lg(y)h(y)dy= \int_{\mathbb{M}}g(y)L^{\ast}h(y)dy,$$
we obtain
\begin{equation*}
\int_{\mathbb{M}}(P_{t}f(y)-f(y))\mu(dy)=\int_{0}^{t}\int_{\mathbb{M}}P_{s}f(y)L^{\ast}\varphi (dy) ds=0
\end{equation*}

Since $\mathcal{C}_{c}^{2}(\mathbb{M})$ is dense in $\mathcal{C}_{0}(\mathbb{M})$ for $\Vert .\Vert_{\infty}$, it follows that $\mu(dy)=\varphi(y)dy$ is an invariant probability as stated.

\underline{Uniqueness of the invariant probability.}
In order to do this, we begin by showing that $\mu$ is an ergodic probability; that is, if a subset $A\subset\mathbb{M}$ satisfies $P_{t}\textbf{1}_{A}=\textbf{1}_{A}$ $\mu -a.s$ for all $t\geqslant 0$, then $\mu(A)$ is either 0 or 1.\\
Let us denote by $f$ the function $P_{t}\textbf{1}_{A}$. Then $f(y)\in\{0,1\}$ for $\mu$-almost $y\in\mathbb{M}$ and $f$ is continuous by point 1 of Theorem \ref{t2}. Since $\mathbb{M}$ is a connected space and $\mu$ has full support, it follows that $f$ is either equal to 0 or 1; and therefore $\mu $ is ergodic.\\
Since two distinct ergodic probabilities are mutual singular, the strong Feller property imply that they must have disjoint support. Since $\mu$ has the whole space, which is connected, as support, the uniqueness of $\mu$ follows.
The second part of the statement is Theorem 4.(i) in \cite{Ichihara}.

The proof that the process is Harris recurrent follows from the proof's lines of Proposition 5.1 in \cite{Ichihara}; which also proves the third point.

\section{Exponential decay in $L^{2}(\mu)$}\label{SectL2cvge}
The goal of this section is to prove the exponential decay in the $L^{2}(\mu)$ norm. The proof heavily relies on the hypocoercitivity method analyzed
 by M.Grothaus and P.Stilgenbauer in \cite{Stilg} whose roots lie in the series of paper \cite{DKMS12}, \cite{DMS10} and \cite{Grothaus}
  initiated by J.Dolbeault, C. Mouhot and C. Schmeiser.
	
	We emphasize that in the particular case where $M=\mathbb{S}^d, n = d+1$ and  $(e_j)_{j=1,\cdots,d+1}$ are the eigenfunctions  associated to the first non-zero eigenvalue, our model coincides with the one studied in section 3 in \cite{Stilg}. 

  For an operator $T$ on some Hilbert space $H$, we denote by $D(T)$ its domain and $T^{\ast}$ its adjoint.
 We begin to recall the \textbf{Data (D)} and  \textbf{Hypotheses (H1)-(H4)} introduced in \cite{Stilg}. For convenience we have chosen to replace certain hypotheses from \cite{Stilg} by slightly stronger ones (see the remark \ref{rem:markovcase} below) which are sufficient for our purpose.
 \begin{defin} (\textbf{The Data (D)})
 \label{hyp:D}
 Let $H$ be a real Hilbert space  and let $(P_t)$ be a strongly continuous semigroup on $H$ with generator $(\mathcal{L},D(\mathcal{L}))$ and core $D\subset D(\mathcal{L}).$  We suppose that
 \begin{enumerate}
 \iti 
 There exist  a closed symmetric  operator $(S,D(S))$ and a closed antisymmetric operator $(A,D(A))$ such that
  $D\subset D(S)\cap D(A),\; A (D)  \subset D$ and $\mathcal{L}_{\vert D}=S_{\vert D}-A_{\vert D}.$
 \itiii 
 There exists a closed subspace $F \subset D(S)$ such that $S_{\vert F} = 0$ and $P(D) \subset D$ where $P$ is the orthogonal projection
 $P : F \oplus F^\perp \rar F: f + g \mapsto f$ for all $(f,g)\in F \times F^\perp$.  
 \end{enumerate}
 \end{defin}


     By density of $D\subset D(A)$, closedness of $A$  and the fact that $P(D) \subset D \subset D(A),$  $AP$ is closed and densely defined. Hence, by Von Neumann's Theorem, $(AP)^\ast AP$ is self-adjoint, closed and densely defined. Thus $(I+(AP)^\ast AP): D((AP)^\ast AP)\rightarrow H$ is invertible with bounded inverse. Set
 \begin{equation} B_0=(I+(AP)^\ast AP)^{-1}(AP)^\ast \text{ on } D((AP)^\ast AP). \end{equation}
 In the following we let $( , )_H$ denote the inner product on $H$ and $\|\cdot\|_H$ the associated norm.
 \begin{defin}(\textbf{Hypotheses (H1)-(H4)})
 \label{hyp:H}
 \begin{enumerate}
 \item[(H1)] $PAP_{\vert D}=0$
 \item[(H2)] (\textsf{Microscopic coercivity}). There exists $\Lambda_{1}>0$ such that for all $ f\in D \cap F^{\perp}$, $$(-Sf,f)_{H}\geqslant \Lambda_{1}\Vert f\Vert_{H}^{2}. $$
 \item[(H3)] (\textsf{Macroscopic coercivity}). There exists $\Lambda_{2}>0$ such that for all $ f\in D((AP)^{\ast}(AP)) \cap F$, \begin{equation}\label{H3}\Vert Af\Vert_{H}^{2}\geqslant \Lambda_{2}\Vert f\Vert_{H}^{2}. \end{equation}

 \item[(H4)] (\textsf{Boundedness of auxiliary operators}). The operators $(B_0S,D)$ and $(B_0A(I-P),D)$ are bounded and there exists constants $N_{1}$ and $N_{2}$ such that for all $f \in F^\perp \cap D$
 \begin{enumerate}
 \item[(H4, a)]
 \begin{equation}\Vert B_0Sf \Vert_{H} \leqslant N_{1}\Vert f \Vert_{H}\end{equation} and
 \item [(H4, b)]
 \begin{equation}\label{H4.2}\Vert B_0 A f\Vert_{H} \leqslant N_{2}\Vert f\Vert_{H} \end{equation}.
 \end{enumerate}
 \end{enumerate}
 \end{defin}
 If furthermore $(I-PA^2P)(D)$ is dense in $H$, then conditions $(H3)$ and $(H4, b)$ are implied by the following conditions, as shown by
 Corollary 2.13 and Proposition 2.15 in \cite{Stilg}.
 \begin{enumerate}
 \item[(H3')] Equation (\ref{H3}) holds for all $f \in D \cap F.$
 \item[(H4')] b) For all $f\in D \cap F$
 \begin{equation} \label{H4.2prim} \Vert A^2 f \Vert_{H}\leqslant N_2  \Vert g \Vert_{H} \end{equation}
 where  $g=(I-PA^2P)f.$
 \end{enumerate}
 \begin{theorem}[Theorem 2 in \cite{DMS10}, Theorem 1 in \cite{DKMS12}, Theorem 2.18 in \cite{Stilg}]\label{theoremrappel} Assume that the assumptions of Definitions \ref{hyp:D} and \ref{hyp:H} hold. Then there exist constants $\kappa_{1},\kappa_{2}\in (0,\infty)$ explicitly computable such that for all $g\in H$ and $t\geqslant 0$,
 \begin{equation} \Vert P_{t}g \Vert_H \leqslant \kappa_{1}e^{-\kappa_{2}t}\Vert g \Vert_H
 \end{equation}
 \end{theorem}
 \begin{rem}
 Following the proof's line of section 3.4 in \cite{DKMS12} and the beginning of the proof of Theorem 2.18 in \cite{Stilg}, one obtains
 \begin{equation}\kappa_1=\sqrt{\frac{1+\varepsilon_{\eta}}{1-\varepsilon_{\eta}}}\leqslant \sqrt{1+ 2\eta} \text{ and } \kappa_2=\varepsilon_{\eta}\frac{\Lambda_{2}}{4(1+\Lambda_{2})},\end{equation}
 with
 \begin{equation} \varepsilon_{\eta}=\frac{\eta}{1+\eta}\frac{\varepsilon_{0}}{\max(1,\varepsilon_{0})},\; \eta >0\end{equation}
 and
 \begin{equation}\label{eps0}
 \varepsilon_{0}=\frac{2\Lambda_{2}\Lambda_{1}}{(1+\Lambda_{2})(2+(1+N_{1}+N_2)^{2})}.
 \end{equation}
 \end{rem}
 \begin{rem}
 \label{rem:markovcase}  In case $(P_t)$ is a Markov semigroup with invariant probability $\mu,$ inducing a strongly continuous semigroup on $L^2(\mu)$,  a natural choice for $H$ is
 $$L^2_0(\mu) = \{f \in L^2(\mu) : \: \int f d\mu = 0\}.$$ This choice will be adopted later. In this case, conditions $(D6)$ and $(D7)$ from  \cite{Stilg} are automatically satisfied and Theorem \ref{theoremrappel} implies that for all $f \in L^2(\mu)$
 $$\Vert P_{t}f - \int f d\mu \Vert_{L^2(\mu)} \leqslant \kappa_{1}e^{-\kappa_{2}t}\Vert f - \int f d \mu \Vert_{L^2(\mu)}.$$
 \end{rem}


 \subsection{\textbf{Application to the Proof of Theorem \ref{t4}}}

 Throughout we
 let $$H=L_0^2(\mu):=\{f\in L^2{(\mathbb{M},\mu)}\: : \int_\mathbb{M} f(y)\mu(dy)=0\}$$ and
  $$L_0^2(e^{-\Phi}) = \{f \in L^2(\RR^n, e^{-\Phi}) \: : \int_{\RR^n} f(u) e^{-\Phi(u)} du = 0\}$$ where $\mu$ and $\Phi$ are like in Definition \ref{probinv}.
  Both $H$ and $L_0^2(e^{-\Phi})$ are equipped with the associated $L^2$ inner product and norm.

  The map $\imath:L_0^2(e^{-\Phi})\hookrightarrow H$   defined by  $\imath(g)(x,u)=g(u)$ injects isometrically $L_0^2(e^{-\Phi})$ into $H.$ We let $$F=\imath(L_0^2(e^{-\Phi}))$$
 and $P :  F \oplus F^{\perp} \rightarrow F$ denote the orthogonal projection onto $F.$
 Alternatively $P$ can be defined as
 \begin{equation}
 \label{eq:defProj} (Pf)(x,u)=\int_M f(x,u)\nu (dx).
 \end{equation}
 Using the notation introduced in section \ref{result} we let  $(P_t)$ denote the semigroup defined by
 $$P_t f(y) = \mathbb{E}(f(Y_t^y))$$ for every bounded  Borel map $f : \mathbb{M} \rightarrow \RR;$ where $(Y_t^y)$ stands for the solution to
 $(\ref{eq5})$ with initial condition $Y_0^y = y.$
 \begin{lem}
  \label{lem:Feller2} $(P_t)$  induces   a strongly continuous contraction semigroup on $H.$
 \end{lem}
 \begin{proof1} By invariance of $\mu$ and Jensen inequality $P_t$ defines a bounded operator on $H$ with norm less than $1$ (as already  proved in Remark \ref{rem5}).

  Let $\varepsilon >0$ and $f\in L^2(\mu)$. By density of $\mathcal{C}_0(\mathbb{M})$ in $L^2(\mu)$, there exists $g\in\mathcal{C}_0 (\mathbb{M})$ such that
 $\Vert f-g\Vert_{L^2(\mu)}<\varepsilon. $
 Thus, by the contraction property
 \begin{eqnarray*}
 \Vert P_t f-f\Vert_{L^{2}(\mu)} &\leqslant &\Vert P_t f-P_t g\Vert_{L^{2}(\mu)}+\Vert P_t g-g\Vert_{L^{2}(\mu)}+\Vert g-f\Vert_{L^{2}(\mu)}\\
 &\leqslant & 2\varepsilon +  \Vert P_t g-g\Vert_\infty.
 \end{eqnarray*}
 Hence, by Feller continuity of $(P_t)$ (see Lemma \ref{l9})
 \begin{equation*}
 \limsup_{t\rightarrow 0} \Vert P_t f-f\Vert_{L^{2}(\mu)}\leqslant 2\varepsilon.
 \end{equation*}
 \end{proof1}
 \begin{rem}
  \label{rem:Feller2}
  Note that the conclusion of Lemma \ref{lem:Feller2} hold true for any Feller Markov semigroup having $\mu$ as invariant measure. This will be used later. \end{rem}
 Let $({\cal L}, D({\cal L}))$ denote the infinitesimal generator of $(P_t)$ (now seen as a strongly continuous semigroup on $H$) and let
  $$D = \mathcal{C}_{c}^{\infty}(\mathbb{M}) \cap H.$$
 \begin{prop}
 \label{prop:data} There exist a closed symmetric operator $(S,D(S))$ and a closed antisymmetric operator  $(A,D(A))$ such that
 \bdes \iti $D$ is a core for $S, A$ and ${\cal L}$ invariant under $S, A, {\cal L}$ and $P.$
 \itii $F \subset D(S)$ and $S|_F = 0.$
 \itiii For all $f \in D$
 \begin{equation}\label{LAPLAC} S(f)=\frac{\sigma^{2}}{2}\Delta_{M}f,\end{equation}
 \begin{equation}\label{A} A(f)  = - G_0(f) = \sum_{i = 1}^n a_{j}u_{j}(\nabla e_{j}(x), \nabla_{x} f)_{TM} -e_{j}(x)\partial_{u_{j}}f\end{equation}
 and \begin{equation}\label{Ldroit} \mathcal{L}f= Lf = Sf-Af.\end{equation}
 \edes
 \end{prop}
 This later proposition shows that conditions of Definition \ref{hyp:D} are fulfilled.

 \medskip
 Let $\eta_1(M) = \eta_1$ denote the {\em spectral gap} of $M.$ That is
 \begin{equation}
 \label{eq:gap}
 \eta_1(M) :=\inf\{\int_{M}\vert \nabla h\vert^{2}\nu(dx) \:
  : h\in H^{1}(M), \: \int_{M}h^{2}\nu(dx) = 1,   \: \int_{M}h\nu(dx) =0 \}
 \end{equation}
 where $\| h\|^{2} =(h,h)_{TM}$ and $(.,.)_{TM}$ is the scalar product on the tangent bundle.
 By a classical result in spectral geometry, compactness of $M$ ensures that $\eta_1 > 0$ and equals the  smallest non zero eigenvalue of $-\Delta_M.$

 \begin{prop}
 \label{prop:Hyp}
  Hypotheses (H1)-(H4) in Definition \ref{hyp:H} hold with
 $$\Lambda_1 = \frac{\eta_1 \sigma^2}{2}, \Lambda_2 = \min_{i = 1, \ldots, n} |\lambda_i| a_i,$$
 $$N_1 = \frac{\sigma^{2}}{2}\sum_{j=1}^{n}\vert\lambda_{j}\vert,$$ and
 $$N_2 = 2 \frac{n}{\min\{\vert \lambda_j\vert , j=1,\cdots,n\}}\sup_{i = 1, \ldots, n} \|\nabla e_i\|_{\infty}^2  \sqrt{4+\sum_{i=1}^n |\lambda_i| a_i} + 4 \|\sum_i e_i^2\|_{\infty}$$
 \end{prop}
\begin{rem} Since $N_1\geqslant \frac{n\sigma^2}{2}\eta_1$, then $2\Lambda_1 < 2+(1+N_{1}+N_2)^{2}$. Hence $\varepsilon_0<1$, where $\varepsilon_0$ is defined by (\ref{eps0}).
\end{rem}
 \subsection{\textbf{Proof of Propositions \ref{prop:data} and \ref{prop:Hyp}}}
 \subsubsection*{\textbf{Proof of Proposition \ref{prop:data}}}
 We first recall some  classical results that will be used throughout.
 \begin{prop}(see e.g~Corollary 1.6, Proposition 2.1, Proposition 3.1, Proposition 3.3 in \cite{EthierKurtz})\label{ResultElem} Let $K$ be the generator of a strongly continuous contracting semi-group $(T_t)_t $ on some Banach space $\mathcal{H}$. Then
 \begin{enumerate}
 \item $K$ is closed and densely defined.
 \item The resolvent set of $K$ contains $(0,\infty)$ and $(\lambda I -K)^{-1}g=\int_0^{\infty}e^{-\lambda t}T_t g dt $, for all $g\in\mathcal{H}$ and $\lambda >0$.
 \item A subspace $D$ of $D(K)$ is a core for $K$ if and only if it is dense in $\mathcal{H}$ and $(\lambda I -K)(D)$ is dense in $\mathcal{H}$ for some $\lambda >0$.
 \item Let $D$ be a dense subset of $\mathcal{H}$ such that $D\subset D(K)$. If $T_t(D)\subset D$ for all $t\geqslant 0$, then $D$ is a core for $K$.
 \end{enumerate}
 \end{prop}
 Similarly to $(P_t),$ let
 $(P_t^S)$ and  $(P_t^A)$
 be the semigroups respectively  induced by  the following stochastic and ordinary differential equation on $\mathbb{M}$:
 $$dY_t^S= \sum_{j=1}^{N}G_{j}(Y_t^S)\circ dB_t^j,$$ and
 \beq \label{edo1} \frac{dY_t^A}{dt} = - G_0(Y_t^A).
 \eeq
  Note that $(P_t^A)$
  is not merely a semigroup but a group of transformation defined as
 \begin{equation}
 P_t^Af(y)=(f\circ \psi_t)(y)
 \end{equation}
 where $\{\psi_t\}$  is the flow induced by (\ref{edo1}). The proofs given in Lemma \ref{l9}, Proposition \ref{p5} and Remark \ref{rem:Feller2} show that, not only $(P_t)$ but also $(P_t^S)$ and $(P_t^A)$ are Feller, leave $\mathcal{C}^2_c(\mathbb{M})$ invariant and admit $\mu$ as invariant probability. Thus, by Remark \ref{rem:Feller2}   and Proposition \ref{ResultElem} they induce strongly continuous semigroups on $H$ whose generators, denoted $S$ and $A$ are closed, densely defined and admit $\mathcal{C}^2_c(\mathbb{M}) \cap H$ as a core.

 Since for all $f \in F, P_t^S f = f,$ assertion $(ii)$ of Proposition \ref{prop:data} is satisfied. Furthermore, the definition of ${\cal L}, A$ and $S$ easily imply assertion $(iii)$ as well as invariance of $D$ under the generators and under $P.$  The end of the proof is given by the two following lemmas.

 \begin{lem} \label{p1}  $D$ is a core for  $\mathcal{L}, S$ and $A.$
 \end{lem}
 \begin{proof1} Let $G$ be one of the operators  ${\cal L}, S$ or $A.$
 It is easily checked that for all $f \in C^2_c(\mathbb{M})$
 $$\|A f\|_{L^2(\mu)} \leq  C \|\nabla f\|_{\infty}$$
 and $$\| S f\|_{L^2(\mu)} \leq \frac{\sigma^2}{2}\|\Delta_M f\|_{\infty}$$ for some $C > 0$ independent of $f.$  Thus $G$ maps continuously  the space
 $\mathcal{C}^2_c(\mathbb{M}) \cap H$ equipped with the $\mathcal{C}^2$ strong topology, into $H$.
 By standards approximation results $\mathcal{C}^\infty_c(\mathbb{M})$ is dense into $\mathcal{C}^2_c(\mathbb{M})$ for the $\mathcal{C}^2$
 strong topology (see e.g~ \cite{Hirsch74}, Chapter 2).
  Since $\mathcal{C}^2_c(\mathbb{M}) \cap H$ is a core for
 $G,$ $(I- G) (\mathcal{C}^2_c(\mathbb{M}) \cap H)$ is dense in $H$ (see Proposition \ref{ResultElem}). Thus $(I-G)(D)$ is dense in $H$ and $D$ is a core.
  \end{proof1}
 \begin{lem}
 \label{symandantisym} $S$ is symmetric and $A^\ast =-A$.
 \end{lem}
 \begin{proof1}
 Let $f,g\in D$. Then
 \begin{eqnarray*}
 (Sf,g)_H &=&\frac{\sigma^{2}}{2}\int\int_{M}(\Delta_{M}f)g \nu(dx) e^{-\Phi}d\xi=-\frac{\sigma^{2}}{2}\int\int_{M}(\nabla f,\nabla g)_{TM} \nu(dx) e^{-\Phi}d\xi\\
 &=&\frac{\sigma^{2}}{2}\int\int_{M}(\Delta_{M}g) f \nu(dx) e^{-\Phi}d\xi = (f,Sg)_H
 \end{eqnarray*}
 Since $D$ is a core for $S$, this proves the symmetry of $S.$\\
  For $f,g\in H$, we obtain from invariance of $\mu$,
 \begin{eqnarray}
 (P_t^Af, g)_H = \int_{\mathbb{M}}(f\circ \psi_t)(y)g(y)\mu(dy)&=&\int_{\mathbb{M}}f(\psi_t (y)) g(\psi_{-t} \circ \psi_t (y))\mu(dy)\\
 &=& \int_{\mathbb{M}}f(y)(g\circ {\psi}_{-t})(y)\mu(dy).
 \end{eqnarray}
 Hence $(P_t^A)^\ast = P^A_{-t}$. In particular, $((P_t^A)^\ast)$ is strongly continuous and admits $-A$ as infinitesimal generator. Now,  when a semigroup and its adjoint are both strongly continuous, the generator of the adjoint equals the  adjoint of the generator. This follows for instance from   Theorem 1.5 in \cite{Phillips} combined with Proposition  \ref{ResultElem} 2. Thus $A^\ast = -A.$
 \end{proof1}
 \subsubsection*{\textbf{Proof of Proposition \ref{prop:Hyp}}}
 For all $f \in D$ let
 \begin{equation}\label{Aj} A_{j}(f)(x,u)= a_{j}u_{j}(\nabla e_{j}(x), \nabla_{x} f)_{TM} -e_{j}(x)\partial_{u_{j}}f.\end{equation}
 so that
 $Af = \sum_{j = 1}^n A_j f.$ Similarly to $A$, $A_j$ enjoys the same properties as $A$. In particular, it leaves $D$ invariant and is antisymmetric:
 $$ (A_j  f, g)_{L^2(\mu)} = -(f, A_j g)_{L^2(\mu)}$$ for all $f,g \in D.$

 Finally, we introduce the following operators
 \begin{eqnarray}
 T&=& (I+(AP)^{\ast}(AP))^{-1} \mbox{ on } H \\
 B_j&=&- T(P A_j) \mbox{ on } D
 \end{eqnarray}
 where $I$ denotes the identity operator. Recall that $B_0$ was introduced to be the operator
 \begin{equation*} B_0=T(AP)^{\ast} \mbox{ on } D ((AP)^* AP).\end{equation*}
 Hypothesis (H1) is immediate because for all
 $f \in D, A_j P f= - e_j(x)  \partial_{u_{j}}(P f)$ and $\int_{M}e_{j}(x)\nu(dx)=0,$  thus
 $P A_j P f  =  0.$

 Hypothesis (H2) follows directly from the variational definition of the spectral gap (\ref{eq:gap}). Indeed for  all $f\in D \cap F^\perp $
 \begin{eqnarray*}
 -(Sf,f)_{L^{2}(\mu)}&=&-\frac{1}{2}\sigma^{2}\int_{\mathbb{R}^{n}}\int_{M}(\Delta_{M}f) f \nu(dx) e^{-\Phi (u)}du\\
 &=&\frac{1}{2}\sigma^{2}\int_{\mathbb{R}^{n}}\int_{M}\vert\nabla_{x}f\vert^{2} \nu(dx) e^{-\Phi (u)}du \geq \frac{\eta_{1}}{2}\sigma^{2} \Vert f \Vert_{L^{2}(\mu)}^{2}.
 \end{eqnarray*}
  For $k = 1, \ldots, n$ let
 $$\alpha_k = |\lambda_k| a_k,$$ so that $$\Phi(u) = \frac{1}{2} \sum_{k = 1}^n \alpha_k u_k^2 + \ln(C(\Phi)).$$
 Let $(P_t^{OU})$ denote the {\em Ornstein-Uhlenbeck} semi-group on $L^2_0(e^{-\Phi})$ defined as
 \beq
 \label{eq:defPtou}
 P_t^{OU} f(u) = \int f( e^{- \mathsf{diag}(\alpha_i) t} u + \mathsf{diag}(\sqrt{1-e^{-2\alpha_i t}}) \xi) e^{-\Phi(\xi)} d\xi
 \eeq
 or, equivalently, $P_t^{OU} f(u) = \mathbb{E}(f(U_t^u))$ where $U_t^u$ is the solution to the linear equation on $\RR^n$
 $$dU_t^i = -\alpha_i U_t^i dt + \sqrt{2} dB^i_t,\; i =  1 \ldots n,$$ with initial condition $U_0^u = u$ and independent Brownian motions $B^1, \ldots, B^n.$

 Let
 $L_{OU}$ denote the generator of $(P_t^{OU})$ on $L^2_0(e^{-\Phi}).$
  The set $$\tilde{D} = {\cal C}_c^{\infty}(\RR^n) \cap L^2_0(e^{-\Phi})$$ is a core\footnote{This is a classical result and can easily be verified as follows. Formula (\ref{eq:defPtou}) shows that the set $C^{\infty}_b(\RR^n)$ of bounded $C^{\infty}$ functions with bounded derivatives is stable under $(P_t^{OU})$; hence a Core by Proposition \ref{ResultElem}. 
  Furthermore for each $f \in C^{\infty}_b(\RR^n)$  it is easy to construct a sequence  $f_n \in C^{\infty}_c(\RR^n)$
  such that $f_n \rightarrow f$ and $L_{OU} f_n \rightarrow L_{OU} f$ in $L^2(e^{-\Phi}).$}
 $L_{OU}$ and for all $f \in \tilde{D} $
 $$L_{OU} f = - \langle \nabla \Phi, \nabla f\rangle + \Delta f.$$

 The next Lemma is similar to Corollary 2.13 and Proposition 3.13 in \cite{Stilg}, 
 \begin{lem} \label{I-PA2P}\bdes
 \iti  For all $f \in F$ $$P A^2 f = \imath \circ  L_{OU} \circ \imath^{-1}(f)$$
 \itii $(I - PA^2P)(D)$ is dense in $H.$
 \itiii (H3) holds with  $\Lambda_2 = \min \{\alpha_k :  \: k =1 \ldots n \}$
 \edes \end{lem}
 \begin{proof1}
 $(i)$ Let $f\in F \cap D$. Then
 \begin{eqnarray}
 \label{eq:A2f}
 \nonumber A^{2} f&=&\sum_{k=1}^{n}A_{k}(Af)=\sum_{k=1}^n A_{k}(\sum_{j=1}^{n}A_{j}f)=\sum_{k=1}^n A_{k}(\sum_{j=1}^{n} e_{j}\partial_{u_{j}} f)\\
 &=& \nonumber \sum_{k=1}^{n}[(\nabla e_{k},\sum_{j=1}^{n}\nabla_x (e_{j}\partial_{u_{j}} f))_{TM}a_{k}u_{k}-\partial_{u_{k}}(\sum_{j=1}^{n} e_{j}\partial_{u_{j}} f)e_{k}]\\
 &=& \sum_{k,j=1}^{n}\partial_{u_{j}}f a_{k}u_{k}(\nabla e_{k},\nabla e_{j})_{TM}-\sum_{k,j=1}^{n}(\partial^2_{u_{j}u_{k}} f)e_{j}e_{k}
 \end{eqnarray}
 Therefore
 \begin{eqnarray} P A^{2} f&=&\sum_{k,j=1}^n\partial_{u_{j}} f a_{k}u_{k}\int_{M}(\nabla e_{k},\nabla e_{j})_{TM}d\nu-\sum_{k,j=1}^{n}(\partial^2_{u_{j}u_{k}} f)\int_{M}e_{j}e_{k}d\nu \notag \\
 &=& \sum_{j=1}^{n}\partial_{u_{j}} f a_{j}u_{j}\vert \lambda_{j}\vert -\sum_{j=1}^{n}(\partial_{u_{j} u_j} f) =
  \sum_{j=1}^{n}\partial_{u_{j}} f \alpha_j u_{j} -\sum_{j=1}^{n}(\partial^2_{u_{j} u_j} f).
 \end{eqnarray}
 This proves the first assertion.

 $(ii)$ $(I-PA^2P)(D \cap F^\perp) = D \cap F^\perp$ is dense in $F^\perp$ because $F^\perp = (I-P)(H), (I-P)(D) \subset D \cap F^\perp$
  and $D$ is dense. Also,
  $(I-PA^2P)(D \cap F) = \imath (I -L_{OU})(\tilde{D})$ is dense in $F$ because, $\tilde{D}$ being a core for $L_{OU}$,
 $(I -L_{OU})(\tilde{D})$ is dense in $L^2_0(e^{-\Phi}).$ This proves $(ii)$.

 $(iii)$  Using antisymmetry of $A,$  assertion $(i)$ and the Poincar\'e inequality for the Gaussian measure $e^{-\Phi(u)} du$ (see e.g~ \cite{ABC+}, chapter 1) we get that for all $f \in F \cap D,$ $$\|A f\|^2_H  = \|A P f\|^2_H  = (- P A^2 P f, f)_H =  ( \imath (f),  L_{OU} \imath (f) )_{L_0^2(e^{-\Phi})}$$
  $$  \geq \min ({\alpha_i}) \|\imath(f)\|_{L^2_0(e^{-\Phi})}
   = \min ({\alpha_i}) \|f\|^2_{H}.$$ This proves (H3'), hence (H3).
 \end{proof1}
 \begin{lem}\label{l4} For $f\in
  D \cap F$, we have $\Vert A f \Vert_{L^{2}(\mu)}^2=\sum_{k=1}^{n}\Vert A_{k} f\Vert_{L^{2}(\mu)}^2=\Vert \nabla f\Vert_{L^2(\mu)}^2$.
 \end{lem}
 \begin{proof1} Let $f\in D \cap F $. Since $ f$ does not depend on the $x$-variable, $A_{j} f=-e_{j}\partial_{u_{j}} f$. The result follows from the fact that the eigenfunctions $(e_j)_{j=1,\cdots ,n}$ are orthonormal in $L^2(M,dx)$.
 \end{proof1}
 The next Lemma is inspired from Lemma 2.4 in \cite{Stilg}
 \begin{lem} \label{l5} For $j=1,\cdots ,n$ and $f\in D$,
 $$\Vert B_{j} f\Vert_{H}\leqslant \frac{1}{2}\Vert (I-P)f\Vert_{H}. $$
 \end{lem}

 \begin{proof1} 
The proof is quite similar to the proof of Lemma 2.4 in \cite{Stilg}. Let $f\in D$ and define $g=B_j f$. Thus $g\in D((AP)^\ast AP)$ and
\begin{equation}
\label{triv} - P A_j  f= g+((AP)^\ast AP)g.
\end{equation}
Because $(I-PA^2P)(D)$ is dense in $H$ (see Lemma \ref{I-PA2P}.(\textbf{ii})), there exists a sequence $(g_n)\subset D$ such that
\begin{equation}\lim_{n\rightarrow \infty} g_n - PA^2Pg_n = g+(AP)^\ast (AP)g.
\end{equation}
Since $P(D),A(D)\subset D$, it follows from Lemma 2.2 in \cite{Stilg} that
\begin{equation}\label{suitecvge} -P A^2 P g_n = ((AP)^\ast (AP))g_n .
\end{equation}
Thus, by continuity of $T$,
\begin{equation} \lim_{n\rightarrow \infty} g_n = g
\end{equation}
and from (\ref{suitecvge})
\begin{equation} \lim_{n\rightarrow \infty} (AP)^\ast (AP)g_n = (AP)^\ast (AP)g.
\end{equation}
Thus, taking the scalar product of (\ref{triv}) with  respect to $g_n$ on both side provides
 $$\lim_{n\rightarrow \infty} - (PA_j f, g_n)_H  - \|g_n\|^2_H - \|AP g_n\|^2_H = 0.$$
Now, using successively  antisymmetry of $A_j,$ Cauchy Schwarz (and Young) inequalities and Lemma \ref{l5},
\begin{eqnarray}
- (PA_j f, g_n)_H = ((I-P) f, A_j P g_n)_H \leq \|(I-P) f\|_H \| A_j P g_n\|_H \\
\leq \frac{1}{4} \|(I-P) f\|^2_H + \|  A_j P g_n\|^2_H \leq  \frac{1}{4} \|(I-P) f\|^2_H + \Vert APg_n\Vert_H^2
\end{eqnarray}
Thus, letting $n$ tends to $\infty$, leads to
\begin{equation} \Vert g\Vert_{H}^2\leqslant \frac{1}{4}\Vert (I-P)f\Vert_{H}^2.
\end{equation}
 \end{proof1}
 \begin{lem}\label{l7} (H4 a) holds with $N_1 = \frac{\sigma^{2}}{2}\sum_{j=1}^{n}\vert\lambda_{j}\vert.$
 \end{lem}
 \begin{proof1} Let $f\in D\cap F^\bot$. Since $\int_{\mathbb{M}}A_j f(y)\mu(dy)=0$, one has
 \begin{eqnarray*}
 -P A f&=& \sum_{j=1}^{n}-P A_{j}f=\sum_{j=1}^{n}P(a_{j}u_{j}(\nabla e_{j},\nabla_{x} f)_{TM}-e_{j}\partial_{u_{j}}f)\\
 &=& \sum_{j=1}^{n}[\int_{M}(\nabla e_{j},\nabla_{x}f)_{TM}a_{j}u_{j}d\nu-\int_{M}e_{j}\partial_{u_{j}}fd\nu].
 \end{eqnarray*}
 Since $S(D)\subset D)$, then
 \begin{eqnarray*}
 -P A Sf=-\frac{\sigma^{2}}{2}\sum_{j=1}^{n}[\int_{M}(\nabla e_{j},\nabla_{x}\Delta_{M}f)_{TM}a_{j}u_{j}d\nu-\int_{M}e_{j}\partial_{u_{j}}\Delta_{M}fd\nu] .
 \end{eqnarray*}
 Because
 \begin{eqnarray*}
 \int_{M}(\nabla e_{j},\nabla_{x}\Delta_{M}f)_{TM}d\nu&=& -\int_{M}\Delta_M e_{j}\Delta_{M }fd\nu=-\lambda_{j}\int_{M}e_{j}\Delta_{M}fd\nu \\
 &=& \lambda_{j}\int_{M}(\nabla e_{j},\nabla_{x}f)_{TM}d\nu \end{eqnarray*}
 and
 $$\int_{M}e_{j}\partial_{u_{j}}\Delta_{M}fd\nu=\int_{M}e_{j}\Delta_{M}\partial_{u_{j}}fd\nu= \int_{M}\Delta_{M}e_{j}\partial_{u_{j}}fd\nu= \lambda_{j}\int_{M}e_{j}\partial_{u_{j}}fd\nu$$
 for all $j=1,\cdots,n$, it follows that
 $$P A Sf=\frac{\sigma^{2}}{2}\sum_{j=1}^{n}\lambda_{j} (P A_{j})f. $$
 By antisymmetry of $A$ (resp. $A_j$)  and Lemma 2.2 in \cite{Stilg}, for all $g$ in $D$, $(AP)^* g = - P A g$ (resp. $(A_j P)^* f = - P A_j f$). Hence
   $$B_{0}Sf = T (AP)^* S f = - T P A S f =  \frac{\sigma^{2}}{2}\sum_{j=1}^{n} \lambda_{j} B_{j}f. $$
 Applying the triangle inequality, one has $$\Vert B_{0}S f\Vert_{L^{2}(\mu)}\leqslant \frac{\sigma^{2}}{2}\sum_{j=1}^{n}\vert\lambda_{j}\vert \Vert B_{j}f\Vert_{L^{2}(\mu)}$$ and
 the result follows from Lemma \ref{l5}.
 \end{proof1}

 The following estimate can be compared with the a priori estimates obtained in \cite{DMS10} and discussed in Appendix $A1$ of \cite{Stilg} (lemmas A3, A4, A5, A7 and Proposition A6) for a more general elliptic equation.  Note, however, that here we provide  an elementary proof allowing  precise estimates by making use of the $\Gamma$ and $\Gamma_2$ operators combined with the specific form of $L_{OU}$.
 \begin{lem} Let $f\in \tilde{D}$ and  \begin{equation} \label{defg} g=(I-L_{OU})f.\end{equation} Then
 \begin{enumerate}
 \item $\Vert \vert Hess(f)\vert_2 \Vert_{L^2(e^{-\Phi})} \leqslant 4 \Vert g\Vert_{L^2(e^{-\Phi})}$
 \item $\Vert \vert \nabla \Phi \vert_2 .\vert \nabla f\vert_2 \Vert_{L^2(e^{-\Phi})} \leqslant 2\sqrt{4+\sum_{i=1}^n\alpha_i} \Vert g\Vert_{L^2(e^{-\Phi})}$,
 \end{enumerate}
 where $\vert .\vert_2$ stands for the usual Euclidean norm and $\vert Hess(f)\vert_2^2 = \sum_{i j} |\partial_{u_i u_j} f|^2.$
 \end{lem}
 \begin{proof1}
 From (\ref{defg}), we have $f=R_1 g$, where $R_1$ is the resolvent operator of $L_{OU}$. Thus $$\Vert f \Vert_{L^2(e^{-\Phi})} \leqslant \Vert g\Vert_{L^2(e^{-\Phi})}$$ and $$\Vert L_{OU}f \Vert_{L^2(e^{-\Phi})} \leqslant 2\Vert g\Vert_{L^2(e^{-\Phi})}.$$
 Let $\Gamma$ be the \textit{``carr\'e du champs''} operator defined by
 \begin{equation}\Gamma (\psi_1,\psi_2)=\frac{1}{2}[L_{OU}(\psi_1\psi_2)-\psi_2 L_{OU}\psi_1 -\psi_1 L_{OU}\psi_2] \end{equation} and \begin{equation}\Gamma_2 (\psi)=\frac{1}{2}\Gamma (\psi,\psi) -\psi L_{OU}\psi).  \end{equation}
 It is known (see for instance Subsection 5.3.1 in \cite{ABC+}) that
 \begin{enumerate}
 \item[(i)]  $\Gamma (f,f)=\vert \nabla f\vert^2_2$ and
 \item[(ii)] $\Gamma_2(f)=\vert Hess(f)\vert^2 _2+\langle \nabla f, Hess(\Phi) \nabla f\rangle\geqslant \vert Hess(f)\vert^2_2$
 \end{enumerate}
 by positive definiteness of $Hess(\Phi)$. Therefore, by invariance and reversibility of $e^{-\Phi(u)}du$,
 \begin{eqnarray} \Vert \vert \nabla f\vert_2 \Vert_{L^2(e^{-\Phi})}^2 &=& \int \Gamma(f,f) e^{-\Phi(u)}du\notag\\
 &=&\int -fL_{OU}f e^{-\Phi(u)}du\notag\\
 &\leqslant & \Vert f\Vert_{L^2(e^{-\Phi})}\Vert L_{OU}f\Vert_{L^2(e^{-\Phi})} \text{ ( by the Cauchy-Schwarz inequality)}\notag\\
 &\leqslant & 2 \Vert g\Vert_{L^2(e^{-\Phi})}^2
 \end{eqnarray}
 and
 \begin{equation} \int \Gamma_2(f) e^{-\Phi(u)}du = \Vert L_{OU}f\Vert_{L^2(e^{-\Phi})}^2\leqslant 4 \Vert g\Vert_{L^2(e^{-\Phi})}^2.
 \end{equation}
 This last inequality implies $(i)$. Set $h=\vert \nabla f\vert_2$ so that $\partial_{u_j}h=\frac{\partial_{u_j}f\partial_{u_j}^2f}{\vert \nabla f\vert_2}$. Following the line of the proof of Lemma A.18 in \cite{Villani} and noting that $\Delta \Phi=\sum_{i=1}^n\alpha_i$, one obtains
 \begin{equation}\int \vert \nabla \Phi\vert^2_2 h^2 e^{-\Phi}du \leqslant \sum_{i=1}^n\alpha_i \int h^2 e^{-\Phi}du +2\sqrt{(\int \vert \nabla \Phi\vert^2_2 h^2 e^{-\Phi}du)(\int \vert \nabla h\vert^2_2 e^{-\Phi}du)}.
 \end{equation}
 Using the Young's inequality $2ab\leqslant \delta^2 a^2+\frac{b^2}{\delta^2}$ with $\delta^2=1/2$, one has
 \begin{equation}\int \vert \nabla \Phi\vert^2_2 h^2 e^{-\Phi}du \leqslant 2\sum_{i=1}^n\alpha_i \int h^2 e^{-\Phi}du + 4\int \vert \nabla h\vert^2_2 e^{-\Phi}du.
 \end{equation}
 Since
 \begin{eqnarray} \vert \nabla h\vert^2_2&=&\sum_{j=1}^n (\frac{\partial_{u_j}f}{\vert \nabla f\vert_2} )^2 (\partial_{u_j}^2 f)^2\notag\\
 &\leqslant&\sum_{j=1}^n(\partial_{u_j}^2 f)^2\notag\\
 &\leqslant&\vert Hess(f)\vert_2^2 ,
 \end{eqnarray}
 we obtain
 \begin{eqnarray} \Vert \vert \nabla \Phi\vert_2 .\vert \nabla f\vert_2 \Vert_{L^2(e^{-\Phi})}^2&\leqslant&2(\sum_{i=1}^n\alpha_i)\int \vert \nabla f\vert^2_2 e^{-\Phi}du +4 \int \vert Hess(f)\vert^2_2 e^{-\Phi}du\notag\\
 &\leqslant&4(\sum_{i=1}^n\alpha_i + 4) \Vert g\Vert_{L^2(e^{-\Phi})}^2 .
 \end{eqnarray}
 \end{proof1}
 \begin{cor} Hypothesis
  (H4') b) holds with $$N_2 = 2 \frac{n}{\min\{\vert \lambda_j\vert , j=1,\ldots,n\}}\sup_{i = 1, \ldots, n} \|\nabla e_i\|_{\infty}^2  \sqrt{4+\sum_{i=1}^n\alpha_i} + 4 \|\sum_i e_i^2\|_{\infty}$$ \end{cor}
  \begin{proof1}
 Let $f \in F \cap D.$ To shorten notation we identify $f$ and $\imath^{-1}(f) \in \tilde{D}.$ Then equation (\ref{eq:A2f}) and Cauchy-Schwarz inequality implies
 \begin{eqnarray*}|A^2 f| &\leqslant& \sum_{j,k=1}^n \vert \partial_{u_j} f\vert |\nabla e_j|_M \frac{\alpha_k}{\vert \lambda_k\vert}u_k|\nabla e_k|_M+\vert \sum_{k,j=1}^{n}(\partial_{u_{j}u_{k}} f)e_{j}e_{k}\vert\\
&\leqslant& (\sum_{j=1}^n\partial_{u_j} f\vert |\nabla e_j|_M)(\sum_{k=1}^n (\alpha_k u_k)|\nabla e_k|_M)\lambda_\ast +\vert \sum_{k,j=1}^{n}(\partial_{u_{j}u_{k}} f)e_{j}e_{k}\vert\\
&\leqslant& n\lambda_\ast\sqrt{\sum_{i = 1}^n (\partial_{u_i} f |\nabla e_i|_M)^2} \sqrt{ \sum_{i = 1}^n  (\alpha_i u_i)^2 |\nabla e_i|_M^2} + |Hess(f)|_2 (\sum_i e_i^2)\\
&\leqslant& n\lambda_\ast\sup_i \|\nabla e_i\|_{\infty}^2 |\nabla f|_2 |\nabla \Phi|_2 + |Hess(f)|_2  \|\sum_i e_i^2\|_{\infty},\end{eqnarray*}
where $\lambda_\ast =\frac{1}{\min\{\vert \lambda_j\vert , j=1,\cdots,n\}}$. The result then follows from the preceding lemma.
  \end{proof1}
\section{Exponential decay in the total variation norm}\label{SectTVcvge}
The idea for proving the exponential decay in total variation consists on translating our problem to a setting for which the arguments used for the exponential decay in $L^{2}(\mu)$ remain valid.

Let $z_{0}\in\mathbb{M}$. Since for all $t>0$, $P_{t}(z_{0},dz)=p_{t}(z_{0},z)dz$ where $p_{t}(z_{0},.)$ is a smooth function and that the invariant probability $\mu$ has a smooth density $\varphi$, one has
\begin{eqnarray*}\Vert P_{t}(z_{0},dz)-\mu(dz)\Vert_{TV}&=&\int_{\mathbb{M}}\vert p_{t}(z_{0},z)-\varphi(z)\vert dz.
\end{eqnarray*}
Because $\varphi>0$, we can define a function $h(t,z_{0},.)$ by
$$h(t,z_{0},z)=\frac{p_{t}(z_{0},z)}{\varphi(z)} $$
By Proposition \ref{p0}, $P_{t}(z_{0},dz)$ has a compact support, ie $p_{t}(z_{0},.)$ has a compact support. Hence so does $h(t,z_{0},.)$. Moreover the smoothness of $\varphi$ and $p_{t}(z_{0},.)$ implies the smoothness of $h(t,z_{0},.)$. Consequently, $h(t,z_{0},.)\in L^{2}(\mathbb{M},\mu)$ and
\begin{eqnarray}\int\vert p_{t}(z_{0},z)-\varphi(z)\vert dz&=&\int\vert h(t,z_{0},z)-1\vert \mu(dz)\notag\\
&\leqslant&(\int (h(t,z_{0},z)-1)^{2} \mu(dz))^{\frac{1}{2}}\notag\\
&=&\Vert h(t,z_{0},.)-1\Vert_{L^{2}(\mu)}.
\end{eqnarray}
Since $\int_{\mathbb{M}}h(t,z_{0},y)\mu(dy)=1$ for all $t$ and $z_{0}$, we have a similar formulation to the one of Theorem \ref{t2}.

So, in order to give the exponential rate of convergence, we will show that $h(t,z_{0},.)$ is solution to the abstract Cauchy problem $\partial_{t}u(t)=\mathcal{L}_{2}u(t)$ in $L^{2}(\mu)$ where $\mathcal{L}_{2}$ is an operator for which the arguments used for $\mathcal{L}$ remain valid.\\
In the following, we denote by $h_t$ (resp. $p_t$) the function $h_t(z_0,.)$ (resp. $p_t(z_0,.)$)

 Since $\partial_{t}p_{t}(z_{0},.)=L^{\ast}(p_{t}(z_{0},.))$ by Theorem 3.(iii) in \cite{Ichihara} (recall that $L^{\ast}$ is defined by (\ref{adjointLEB})), then
\begin{eqnarray} \partial_{t}h_{t}&=&\frac{\partial_{t}p_{t}}{\varphi}\notag\\
&=&\frac{L^{\ast}(p_{t})}{\varphi}\notag\\
&=&\frac{\sigma^{2}}{2}\Delta_{M}h_{t}+\sum_{k=1}^{n}a_{k}u_{k}\frac{div_{x}(\nabla e_{k}(x) p_{t})}{\varphi}-\sum_{k=1}^{n}\frac{\partial_{u_{k}}p_{t}}{\varphi}e_{k}(x)
\end{eqnarray}
Because $\partial_{u_{k}}\varphi=-a_{k}u_{k}\vert\lambda_{k}\vert \varphi$, 
 $$-\frac{\partial_{u_{k}}p_{t}}{\varphi}=-\partial_{u_{k}}h_{t}+a_{k}u_{k}\vert\lambda_{k}\vert h_{t}.$$
Moreover, $$\frac{div_{x}(\nabla e_{k}(x) p_{t})}{\varphi}=\Delta_{M}(e_{k}) h_{t}+(\nabla e_{k}(x),\nabla_{x}h_{t})_{TM}. $$
Hence,
\begin{eqnarray*} \partial_{t}h_{t}&=&\frac{\sigma^{2}}{2}\Delta_{M}h_{t}+\sum_{k=1}^{n}a_{k}u_{k}(\nabla e_{k}(x),\nabla_{x}h_{t})_{TM} -\sum_{k=1}^{n}\partial_{u_{k}}h_{t} e_{k}(x)\\
&=:&L_{2}h_{t}.
\end{eqnarray*}
Thus, $h_t= T(t-1)h_1$, where $T(t)$ is the semi-group whose infinitesimal generator restricted to $C_c^{\infty}(\mathbb{M})$ is $L_2$.
Because
 $$L_{2}=S+\sum_{k=1}^{n}A_{k},$$
whereas  $$L=S-\sum_{k=1}^{n}A_{k},$$
$L_{2}$ is the adjoint operator of $L$ in $L^{2}(\mu)$.
So all the arguments used for proving Theorem \ref{t4} for $L$ work for $L_{2}$.\\
Applying Theorem \ref{t4} to $L_{2}$ with $g_{t}=h_{t+1}$ gives the result.
\appendix
\section{A deterministic study}
In this Appendix, we study on $\mathbb{S}^{1}\times\mathbb{R}^{2}$ the ODE
\begin{equation}\label{eqA1}
\left\{
\begin{array}{r c l}
\dot{X}_{t}&=& (\sin(X_{t})U_{t}-\cos(X_{t})V_{t})\\
\dot{U}_{t}&=& \cos(X_{t}) \\
\dot{V}_{t}&=& \sin(X_{t})
\end{array}
\right.
\end{equation}
in order to prove Theorem \ref{theoremCercle}.
 Since the vectorial field $F$ defined by
\begin{equation}
F(X,U,V)=\begin{pmatrix}
(\sin(X)U-\cos(X)V)\\
\cos(X)\\
\sin(X)
\end{pmatrix}
\end{equation}
is smooth and sub-linear,it induces a smooth flow $\psi:\mathbb{R}\times (\mathbb{S}^{1}\times\mathbb{R}^{2})\rightarrow \mathbb{S}^{1}\times\mathbb{R}^{2}$.
A first and important observation is
\begin{prop}\label{pA0} If the initial condition for the ODE (\ref{eqA1}) is
$$(X_{0},U_{0},V_{0})=(X_{0},\cos(X_{0}),\sin(X_{0})),$$ then $$\psi_{t}(X_{0},U_{0},V_{0})=(X_{0},\cos(X_{0})(t+1),\sin(X_{0})(t+1)) \: \forall t\in\mathbb{R}.$$
 In particular, the line
$$\{(X,Y,Z)\in\mathbb{S}^{1}\times\mathbb{R}^{2}\: : \: X=X_{0},\:\exists t\in\mathbb{R}\: \text{ such that } (Y,Z)=(\cos(X_{0})t,\sin(X_{0})t)\}$$
 is invariant under $\psi$.
\end{prop}
\begin{proof1} By the hypothesis, we have $\dot{X}(0)=0$. Hence $X(t)=X_{0}$ for all $t\in\mathbb{R}$. Therefore, $U(t)=\cos(X_{0})(t+1)$ and $V(t)=\sin(X_{0})(t+1)$
\end{proof1}
An immediate consequence is
\begin{cor} \label{corA1} If $\dot{X}(0)>0$ (respectively $\dot{X}(0)<0$), then $\dot{X}(t)>0$ (respectively $\dot{X}(t)<0$) for all $t$.
\end{cor}
\begin{proof1} We proceed by contradiction. Hence, by continuity of $\dot{X}$, there exists $t_{0}$ such that $\dot{X}(t_{0})=0$. Then the two last Propositions imply that $\dot{X}(t)=0$ for all $t$. In particular $\dot{X}(0)=0$, which is a contradiction.
\end{proof1}
Let
\begin{equation}
\begin{pmatrix}
x\\
u\\
v \end{pmatrix}=
\Xi ( \begin{pmatrix}
X\\
U\\
V \end{pmatrix})=  \begin{pmatrix}
X\\
\cos(X)U+\sin(X)V\\
-\sin(X)U+\cos(X)V \end{pmatrix}.
\end{equation}
Note that $(u,v)$ is obtained from $(U,V)$ by a rotation of angle $-X$. Then, in the new variable,
the ODE (\ref{eqA1}) becomes the ODE
\begin{equation} \label{eqA2}
\dot{x}(t)=- v(t)
\end{equation}
\begin{equation}\label{eqA3}
\left\{
\begin{array}{r c l}
\dot{u}(t)&=&1- v(t)^{2} \\
\dot{v}(t)&=&  u(t)v(t)
\end{array}
\right.
\end{equation}
Let \begin{equation}
H(u,v) = \left\{ \begin{array}{l}
             \frac{1}{2}(u^2 + v^2 - \log(v^2)),  \mbox{ if } v \neq 0,\\
             \infty, \mbox{ if } v = 0.
           \end{array}  \right.
\end{equation}
\begin{prop} The function $H$ is a first integral for the ODE (\ref{eqA3}).
\end{prop}
\begin{proof1}
Let $v_{0}\neq 0$. Deriving $H$ with respect to $t$ and applying the chain rule, we obtain
\begin{eqnarray*}
\frac{d}{dt} H(u,v)&=& (u\dot{u}+v\dot{v})-\frac{\dot{v}}{v}\\
				&=& (u-uv^{2}-vuv)-u\\
				&=& 0
\end{eqnarray*}
\end{proof1}

Note that $H$ is convex, reaches its global minimum in $(0,\pm 1)$ and takes the value $1/2$ at these points.\\
For $c\in [1/2,\infty[$, let $$H_{c}^{+}=H^{-1}(c)\cap \{ v>0\} ,\; H_{c}^{-}=H^{-1}(c)\cap\{v<0\}$$ and set $H_{\infty}=\{v=0\}$. Then, we define $\mathbb{T}_{c}^{\alpha}=\mathbb{S}^{1}\times H_{c}^{\alpha}$ for $\alpha\in \{+,-\}$ and $T_{\infty}=\mathbb{S}^{1}\times H_{\infty}$.

Since the function $H$ is strictly convex on $\{v>0\} $ and $\{ v<0\}$, we observe that $T_{1/2}^{\alpha}$ is a closed curve, $T_{c}^{\alpha}$ a torus and $T_{\infty}$ a cylinder.

A first result is
\begin{prop}\label{pA1} Let $(x(t),u(t),v(t))$ be a solution of the ODE defined by (\ref{eqA2}) and (\ref{eqA3}).
\begin{enumerate}
\item[(i)] $\mathbb{T}_{1/2}^{\alpha}$ is a periodic orbit with period $2\pi$, $\alpha\in\{+,-\}$
\item[(ii)] On $T_{\infty}$, the dynamic takes the form $(x(t),u(t),v(t))=(x(0),u(0)+t,0).$
\end{enumerate}
For $c>1/2$, let $T_{c}$ be the period of (\ref{eqA3}) on $H_{c}^{\alpha}$
\begin{enumerate}
\item[(iii)] If $\frac{x(T_{c})}{2\pi}\in\mathbb{Q}$, then every trajectory on $T_{c}^{\alpha}$ is periodic with period $qT_{c}$ if the irreducible fraction of $\frac{x(T_{c})}{2\pi}$ writes $\frac{p}{q}$.
\item[(iv)] If $\frac{x(T_{c})}{2\pi}\notin\mathbb{Q}$, then every trajectory on $\mathbb{S}^{1}\times H^{-1}(c)$ is dense either on $T_{c}^{+}$ or $T_{c}^{-}$.
\end{enumerate}
\end{prop}
\begin{proof1}
Points $(i)$ and $(ii)$ follow immediately from (\ref{eqA2}), (\ref{eqA3}) and the function $H$.

Without loss of generality, we assume that $x(0)=0$. Let $c>1/2$. Because for $m\in\mathbb{N}^{\ast}$, we have
\begin{eqnarray}\label{periodicity}
 x(mT_{c})&=&\int_{0}^{mT_{c}}\dot{x}(t)dt=-\int_{0}^{mT_{c}}v(t)dt\notag\\
 &=&-m\int_{0}^{T_{c}}v(t)dt\notag\\
 &=&m\int_{0}^{T_{c}}\dot{x}(t)dt,\notag\\
 &=&m x(T_{c})
\end{eqnarray}
we obtain that when $(u(t),v(t))$ is back to its initial condition, then $x(t)$ does a rotation of angle $x(T_{c})$. Hence if $\frac{x(T_{c})}{2\pi}=\frac{p}{q}$, with $q\in\mathbb{N}^{\ast}$, $p\in\mathbb{Z}$ and such that the fraction is irreducible, then
\begin{eqnarray*}
 2p\pi&=& qx(T_{c})\\
 &=& x(qT_{c}).
\end{eqnarray*}
This proves $(iii).$

If $\frac{x(T_{c})}{2\pi}\notin\mathbb{Q}$, then $(x(qT_{c}))_{q\in\mathbb{N}}$ is dense on $\mathbb{S}^{1}$. Now, assume without lost of generality that $v(0)<0$ and let $T$ be the first time such that $x(T)=2\pi$. We claim that $(u(nT),v(nT))_{n\in\mathbb{N}}$ is dense on $ H_{c}^{-}$. Indeed, if it is not the case, then it is periodic since $ H_{c}^{-}$ is a closed simple curve. This implies that $(x(t),u(t),v(t))$ is periodic with period $n_{0}T$. Thus, there exists $q\in\mathbb{N}$ such that $n_{0}T=qT_{c}$. Therefore, by (\ref{periodicity}), we have $2n_{0}\pi=x(qT_{c})=qx(T_{c})$; so that $\frac{x(T_{c})}{2\pi}=\frac{n_{0}}{q}$. This is a contradiction.

The density of $(x(qT_{c}))_{q\in\mathbb{N}}$ on $\mathbb{S}^{1}$ and the one of $(u(nT),v(nT))_{n\in\mathbb{N}}$ on $ H_{c}^{-}$ implies the density of $((x(t),u(t),v(t)))_{t\geqslant 0}$ on $T_{c}^{-}$. This proves $(iv)$.
\end{proof1}
From now, we assume without lost of generality that $v(0)<0$ (the case $v(0)>0$ being symmetric).
In order to derive properties of $c\mapsto T_{c}$ (see Proposition (\ref{pA1})), we change the time scale by use of $t\mapsto x(t)$. This is possible because it is strictly increasing. We denote by $y$ the inverse function of $x$. Since we have assumed that $x(0)=0$, it follows that $y(0)=0$.\\
Set $u_{2}(t)=u(y(t))$ and $v_{2}(t)=v(y(t))$. Therefore $(u_{2},v_{2})$ is solution to the ODE
\begin{equation}\label{eqA5}
\left\{
\begin{array}{r c l}
\dot{u}_{2}(t)&=&( v_{2}(t)-\frac{1}{v_{2}(t)}) \\
\dot{v}_{2}(t)&=&- u_{2}(t)
\end{array}
\right.
\end{equation}
with initial condition $(u(0),v(0))$. Observe that $H$ is still a first integral for this system.
\begin{prop}\label{pA2} Let $(x(t),u(t),v(t))$ be a solution to the ODE defined by equation (\ref{eqA2}) with initial condition $(0,u_{0},v_{0})$ and let $(t,u_{2}(t),v_{2}(t))$ where $(u_{2}(t),v_{2}(t))$ is the solution to the ODE defined by equation (\ref{eqA5}) with initial condition $(u_{0},v_{0})$. \\
Then $(x(t),u(t),v(t))$ is periodic in $\mathbb{S}^{1}\times\mathbb{R}^{2}$ iff $(t,u_{2}(t),v_{2}(t))$ is periodic in $\mathbb{S}^{1}\times\mathbb{R}^{2}$.\\
 Further, if $T$ is the period of $(x(t),u(t),v(t))$, then $x(T)$ is the period of $(t,u_{2}(t),v_{2}(t))$.
\end{prop}
\begin{proof1} Straightforward.
\end{proof1}
Denote by $T_{c,2}$ the period of $(u_{2}(t),v_{2}(t))$, where $c=H(u_{2}(0),v_{2}(0))>1/2$. Then\begin{equation} T_{c,2}= x(T_{c}).\end{equation} \\
An immediate consequence of Propositions \ref{pA1} and \ref{pA2} is that $(t,u_{2}(t),v_{2}(t)) $ is periodic if and only if
\begin{equation}\label{critere} \frac{T_{c,2}}{2\pi}\in\mathbb{Q}. \end{equation}

In the rest of this Appendix, we study the \textit{"period-function"} \begin{equation}\label{perfct} f:(1/2, +\infty)\rightarrow \mathbb{R}_{+}:c\mapsto T_{c,2}.\end{equation}
  First notice that $(0,1)$ and $(0,-1)$ are stationary points for the ODE (\ref{eqA5}).\\
Let $(u_{0},v_{0})\in\mathbb{R}\times (0,\infty)$. By symmetry of $H$ along the line $v_{2} = 0$, what follow remains true for $v_{0}<0$.

 Set $c=H(u_{0},v_{0})$. Since $H$ is a first integral, then $H(u_{2}(t),v_{2}(t))=c$ for all $t$.\\
  Using the fact that $\dot{v}_{2}=-u_{2}$, we have that
\begin{equation}
\frac{1}{2}\dot{v}_{2}^{2}+(\frac{v_{2}^{2}}{2}-\log(v_{2}))=c.
\end{equation}
Set $\phi(v)=(\frac{v^{2}}{2}-\log(v)).$
\begin{center}
\includegraphics[width=0.35\textwidth]{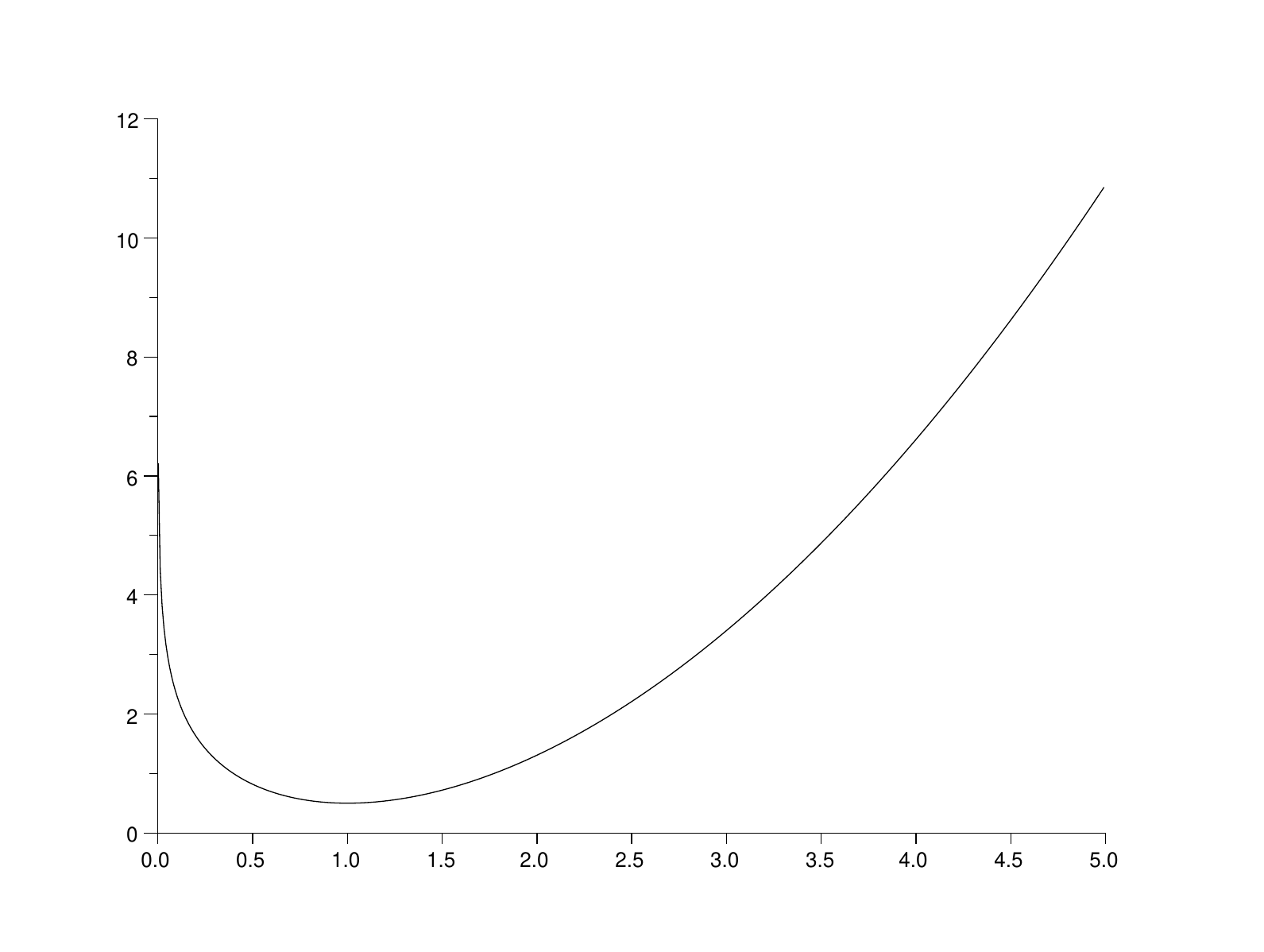}
\captionof{figure}{\scriptsize{Graph of the function $\phi$.}}
\end{center}
Since the curve $H^{-1}(c)$ is symmetric along the line $u_{2} = 0$, we have that
\begin{equation}
\frac{T_{c,2}}{2}=\int_{c_{1}}^{c_{2}}\frac{dv}{\sqrt{2(c-\phi(v))}},
\end{equation}
ie
\begin{equation}
T_{c,2}=\sqrt{2}\int_{c_{1}}^{c_{2}}\frac{dv}{\sqrt{(c-\phi(v))}},
\end{equation}
where $0<c_{1}<1<c_{2}<\infty$ are the roots of the function $v\mapsto\phi(v)-c$.\\
Denote by $h$ the inverse function of $\phi$ restricted to $[1,\infty)$ and by $g$ the inverse function of $\phi$ restricted to $(0,1)$. By a change of variable, we then obtain
\begin{equation}
\int_{1}^{c_{2}}\frac{dv}{\sqrt{(c-\phi(v))}}=\int_{\frac{1}{2}}^{c}\frac{h'(v)dv}{\sqrt{(c-v)}}
\end{equation}
and
\begin{equation}
\int_{c_{1}}^{1}\frac{dv}{\sqrt{(c-\phi(v))}}=-\int_{\frac{1}{2}}^{c}\frac{g'(v)dv}{\sqrt{(c-v)}}.
\end{equation}
Therefore
\begin{equation}\label{conti}
f(c)=T_{c,2}=\sqrt{2}\int_{\frac{1}{2}}^{c}\frac{(h'-g')(v)}{\sqrt{(c-v)}}dv=\int_{\mathbb{R}}\Lambda(v)A(c-v)dv=(\Lambda\ast A)(c),
\end{equation}
where $\ast$ stands for the convolution product, $\Lambda(v)=\sqrt{2}(h'-g')(v)\mathbf{1}_{v>1/2}$ and $A(v)=\frac{1}{\sqrt{v}}\mathbf{1}_{v>0}$.\\
Hence
\begin{equation}
f'(c)=(\Lambda\ast A')(c).
\end{equation}
Since $g(v)\in (0,1)$ and $h(v)>1 $ for $v\in (1/2,c)$, then $g'(v)=\frac{1}{\phi '(g(v))}<0$ and $h'(v)=\frac{1}{\phi '(h(v))}>0$. Using the fact that $A'(v)=-\frac{1}{2}\mathbf{1}_{v>0}\frac{1}{\sqrt{v^{3}}}$, we have
\begin{equation}\label{eqA6}
 f'(c)<0 \text{ for all } 1/2<c<\infty.
\end{equation}
Our next goal is now to study the limiting behaviour $c\rightarrow 1/2$ and $c\rightarrow \infty$
\begin{lem} \label{lA1} Let $c>1/2$ and let $c_{1}$ and $c_{2}$ the two roots of the function $v\mapsto\phi(v)-c$. Then
$$T_{c,2}\geqslant 2\sqrt{2}[\sqrt{\frac{c_{1}}{1+c_{1}}}+\sqrt{\frac{c_{2}}{1+c_{2}}}]. $$
\end{lem}
\begin{proof1} By convexity of $\phi$, we have $\frac{\phi (v)-\phi (c_{1})}{v-c_{1}}\geqslant \phi '(c_{1})$. Hence $$\sqrt{c-\phi(v)}\leqslant \sqrt{-\phi '(c_{1})}\sqrt{v-c_{1}} .$$
Therefore $$\int_{c_{1}}^{1}\frac{dv}{\sqrt{c-\phi(v)}}\geqslant \frac{1}{\sqrt{-\phi '(c_{1})}}\int_{c_{1}}^{1}\frac{dv}{\sqrt{v-c_{1}}}=2\frac{\sqrt{1-c_{1}}}{\sqrt{-\phi '(c_{1})}}. $$
Since $-\phi '(v)=\frac{1}{v}- v$, $-\phi '(c_{1})=(1- c_{1}^{2})/c_{1}$ and thus
$$\int_{c_{1}}^{1}\frac{dv}{\sqrt{c-\phi(v)}}\geqslant 2\sqrt{\frac{c_{1}}{(1+c_{1})}}. $$
Once again convexity of $\phi$ implies $\frac{\phi (c_{2})-\phi(v)}{c_{2}-v}\leqslant \phi '(c_{2})$, so that $c-\phi(v)\leqslant \phi'(c_{2})(c_{2}-v). $
 By proceeding as above, we obtain $$\int_{1}^{c_{2}}\frac{dv}{\sqrt{c-\phi(v)}}\geqslant 2\sqrt{\frac{c_{2}}{(1+c_{2})}} .$$
Hence $$f(c)=T_{c,2}=\sqrt{2}[\int_{c_{1}}^{1}\frac{dv}{\sqrt{c-\phi(v)}}+\int_{1}^{c_{2}}\frac{dv}{\sqrt{c-\phi(v)}}]\geqslant 2\sqrt{2}[\sqrt{\frac{c_{1}}{1+c_{1}}}+\sqrt{\frac{c_{2}}{1+c_{2}}}]. $$
\end{proof1}
\begin{lem} $\lim_{c\rightarrow 1/2} f(c)=\sqrt{2}\pi.$
\end{lem}
\begin{proof1}
We have $c_{1},c_{2}\rightarrow 1$ as $c\rightarrow 1/2$. Thus, it implies that $\log(v)\approx (v-1)-\frac{1}{2}(v-1)^{2}$ for $v\in(c_{1},c_{2})$ and therefore $$\phi(v)=\frac{1}{2}(v-1+1)^{2}-\log(v)\approx \frac{1}{2}+(v-1)^{2}.$$
But
\begin{eqnarray*}
\int_{c_{1}}^{c_{2}}\frac{dv}{\sqrt{c-\frac{1}{2}-(v-1)^{2}}}&=&\frac{1}{\sqrt{c-\frac{1}{2}}}\int_{c_{1}-1}^{c_{2}-1}\frac{dv}{\sqrt{1-(v/\sqrt{c-\frac{1}{2}})^{2}}}\\
&=&\int_{\frac{c_{1}-1}{\sqrt{c-1/2}}}^{\frac{c_{2}-1}{\sqrt{c-1/2}}}\frac{du}{\sqrt{1-u^{2}}}\\
&=& \arcsin (\frac{c_{2}-1}{\sqrt{c-1/2}})+\arcsin(\frac{1-c_{2}}{\sqrt{c-1/2}})
\end{eqnarray*}
Since for $c$ sufficiently close to $1/2$, $c=\phi(1+c_{j}-1)\approx \frac{1}{2}+ (c_{j}-1)^{2}$, then $\lim_{c\rightarrow 1/2}\frac{\vert c_{j}-1\vert}{\sqrt{c-\frac{1}{2}}}= 1, j=1,2$.  \\
Thus, $\lim_{c\rightarrow 1/2}\int_{c_{1}}^{c_{2}}\frac{dv}{\sqrt{c-(v-1)^{2}}}= \pi$ and therefore
\begin{equation}
\lim_{c\rightarrow 1/2}f(c)=\lim_{c\rightarrow 1/2}T_{c,2}=\lim_{c\rightarrow 1/2} \sqrt{2}\int_{c_{1}}^{c_{2}}\frac{dv}{\sqrt{c-\frac{1}{2}-(v-1)^{2}}}= \sqrt{2}\pi .
\end{equation}
\end{proof1}
\begin{rem} One can prove that $\sqrt{2}\pi $ is the period of the orbits from the linear ODE 
\begin{equation}
\left\{
\begin{array}{r c l}
\dot{u}(t)&=& 2v(t)\\
\dot{v}(t)&=&-u(t).
\end{array}
\right.
\end{equation}
But this is nothing else than the linearized system at $(0,1)$ from the ODE (\ref{eqA5}).
\end{rem}
Summarizing all these information concerning $T_{c,2}$, we obtain
\begin{prop}\label{propSUM} The \textit{"period-function"} $f:(1/2,\infty)\rightarrow \mathbb{R}_{+}:c\mapsto T_{c,2}$ is continuous, decreasing, bounded from below by $2\sqrt{2}$ and converge to $\sqrt{2}\pi$ when $c$ tends to $1/2$.
\end{prop}
\begin{proof1} The decreasing property comes from (\ref{eqA6}) whereas the continuity follows from (\ref{conti}). While $c_{1}$ converges to $0$ and $\frac{c_{2}}{1+c_{2}}$ converges to $1$ when $c$ tends to $\infty$, then Lemma \ref{lA1} combined with the decreasing property implies that $f(c)\geqslant 2\sqrt{2}$ for all $c>1/2$.

Since $f$ is decreasing, then $\sup_{c>1/2}f(c)=\lim_{c\rightarrow 1/2}f(c)=\sqrt{2}\pi$.
\end{proof1}
\begin{center}
\includegraphics[width=0.35\textwidth]{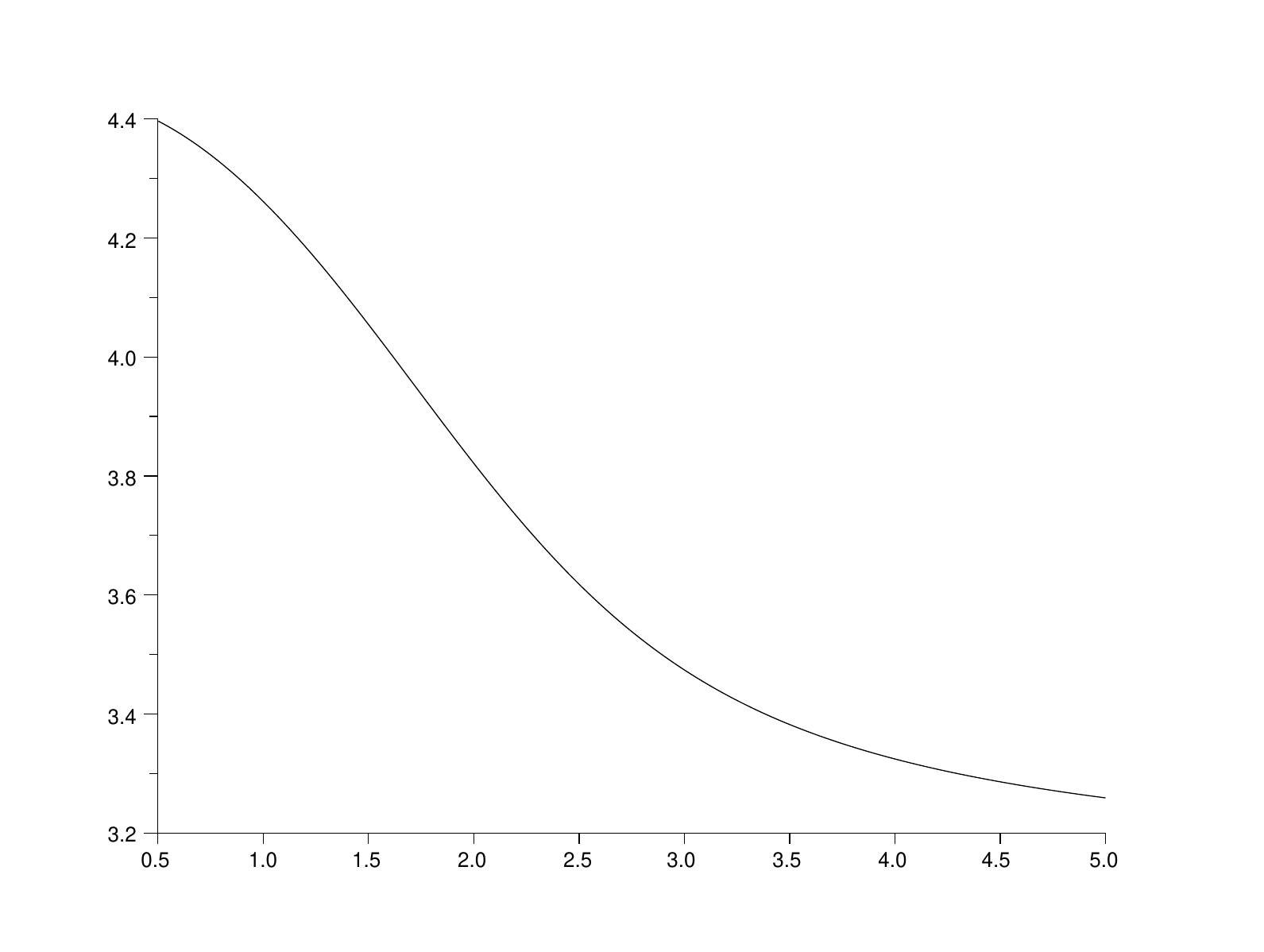}
\captionof{figure}{\scriptsize{Graph of the function $c\mapsto T_{c,2}$.}}
\end{center}

\end{document}